# DECAY ESTIMATES FOR 1-D PARABOLIC PDES WITH BOUNDARY DISTURBANCES

**Iasson Karafyllis[*] and Miroslav Krstic[**]**

[*]Dept. of Mathematics, National Technical University of Athens,
Zografou Campus, 15780, Athens, Greece, email: iasonkar@central.ntua.gr

[**]Dept. of Mechanical and Aerospace Eng., University of California, San Diego, La Jolla, CA 92093-0411, U.S.A., email: krstic@ucsd.edu

**Abstract**

In this work decay estimates are derived for the solutions of 1-D linear parabolic PDEs with disturbances at both boundaries and distributed disturbances. The decay estimates are given in the $L^2$ and $H^1$ norms of the solution and discontinuous disturbances are allowed. Although an eigenfunction expansion for the solution is exploited for the proof of the decay estimates, the estimates do not require knowledge of the eigenvalues and the eigenfunctions of the corresponding Sturm-Liouville operator. Examples show that the obtained results can be applied for the stability analysis of parabolic PDEs with nonlocal terms.

## 1. Introduction

The derivation of decay estimates for the solution of parabolic Partial Differential Equations (PDEs) is a challenging topic, which has attracted the interest of many researchers (see [10,14,15,36,39,40,41,44,48]). The main tool for the derivation of decay estimates is the combination of maximum principles and the so-called "energy" method, i.e., the use of an appropriate functional, which satisfies certain differential inequalities that allow the estimation of the decay rate of the solution. Usually, decay estimates are obtained for systems which do not include time-varying disturbances in the PDE problem.

Recently, the derivation of decay estimates for parabolic PDEs with disturbances was studied by many researchers working mostly in mathematical control theory. Decay estimates for systems with disturbances are related to the Input-to-State Stability (ISS) property (first developed by E. D. Sontag in [47] for systems described by Ordinary Differential Equations-ODEs). The intense interest of researchers in control theory in ISS is justified because: (a) control systems are systems with inputs, and (b) because ISS can be used for the stability analysis by means of small-gain theorems (see Chapter 5 in [19] and references therein). The extension of ISS to systems described by PDEs required novel mathematical tools and approaches (see for example [1,4,5,6,7,16,17,18,20,26,29,30,31,32,43]).

In particular, for PDE systems there are two qualitatively distinct locations where a disturbance can appear: the domain (a distributed disturbance appearing in the PDE) and the boundary (a disturbance that appears in the Boundary Conditions-BCs). Most of the existing results in the literature deal with distributed disturbances. Boundary disturbances present a major challenge, because the transformation of the boundary disturbance to a distributed disturbance leads to decay estimates involving the boundary disturbance and some of its time derivatives (see for example [1]). This is explained by the use of unbounded operators for the expression of the effect of the boundary disturbance (see the relevant discussion in [29] for inputs in infinite-dimensional systems that are expressed by means of unbounded linear operators). Moreover, although the construction of Lyapunov and "energy" functionals for PDEs has progressed significantly during the last years (see for example [2,27,28,30,33,39,40,41,43,44]), none of the proposed Lyapunov functionals can be used for the derivation of the ISS property w.r.t. boundary disturbances. Therefore, the "energy" method cannot provide decay estimates for parabolic PDEs with boundary disturbances. On the other hand, it should be noted that Reaction-Diffusion PDEs with boundary disturbances arise naturally when studying heat and mass transfer phenomena, where flux disturbances appear at the boundaries and the reaction terms are the result of chemical reactions. Parabolic PDEs with boundary disturbances appear also in fluid dynamics (e.g., Navier-Stokes) where boundary/wall disturbances occur naturally in various flow problems.

The recent articles [21,22,23] suggested methodologies for the derivation of decay estimates for 1-D parabolic PDEs with boundary and domain disturbances. Two different methodologies were used in [21,22,23]: the eigenfunction expansion of the solution and the approximation of the solution by means of finite-difference schemes. The obtained ISS estimates were expressed in weighted $L^1, L^2$ and $L^\infty$ norms for the solution under strict regularity requirements for the disturbances ($C^2$ regularity for boundary disturbances and $C^1$ regularity for distributed disturbances) and it was shown that such estimates can be used in a straightforward way for the derivation of decay estimates for parabolic PDEs with nonlocal terms. The interest for the stability analysis of parabolic PDEs with nonlocal terms is strong, both from the PDE literature as well as from the numerical analysis literature (see [8,9,11,12,13,25,37,38]). However, recent advances in feedback control of PDEs has forced the control literature to deal with PDEs containing nonlocal terms. This happened because the feedback law itself is a functional of the solution of the PDE and appears as a nonlocal term either in the BC or in the PDE (see [24,45,46] and references therein).

This paper focuses on 1-D parabolic PDEs with disturbances acting on both boundary sides and distributed disturbances. The contribution of the paper is threefold:
- the derivation of decay estimates in the $L^2$ norm for discontinuous disturbances,
- the derivation of decay estimates in the $H^1$ norm for certain cases, and
- the application of the obtained decay estimates to the stability analysis of parabolic PDEs with nonlocal terms.

More specifically, our first main result (Theorem 2.4) extends recent results (in [21,22,23]) to various directions: discontinuous boundary and domain disturbances can be handled and the obtained decay estimate is less conservative from the corresponding estimates in [21,22,23]. The derivation of decay estimates in the $H^1$ norm (Theorem 2.8 and Theorem 2.9) is achieved for two cases: (a) the case of Dirichlet BCs at both end, and (b) the case of Dirichlet BC on the one end and Robin (or Neumann) BC on the other end. The obtained decay estimates involve the estimation of the principal eigenvalue of a Sturm-Liouville (SL) operator. To this purpose, we develop tools which allow the estimation of the principal eigenvalue (Proposition 2.6).

The structure of the paper is as follows. Section 2 is devoted to the presentation of the problem and the statement of the main results which allow the derivation of decay estimates in various norms (Theorem 2.4, Theorem 2.8 and Theorem 2.9). The application of the obtained decay estimates to the stability analysis of parabolic PDEs with nonlocal terms is illustrated in Section 3. The examples show the exploitation of the decay estimates for the derivation of small-gain conditions for global exponential stability of the zero solution. Section 4 of the present work contains the proofs of all (main and auxiliary) results. The conclusions of the paper are provided in Section 5.

**Notation.** Throughout this paper, we adopt the following notation.
* $\Re_+ := [0, +\infty)$.
* Let $U \subseteq \Re^n$ be a set with non-empty interior and let $\Omega \subseteq \Re$ be a set. By $C^0(U)$ (or $C^0(U;\Omega)$), we denote the class of continuous mappings on $U$ (which take values in $\Omega$). By $C^k(U)$ (or $C^k(U;\Omega)$), where $k \geq 1$, we denote the class of continuous functions on $U$, which have continuous derivatives of order $k$ on $U$ (and also take values in $\Omega$).

* Let $I \subseteq \Re$ be an interval. A function $f : I \to \Re$ is called right continuous on $I$ if for every $t \in I$ and $\varepsilon > 0$ there exists $\delta(\varepsilon, t) > 0$ such that for all $\tau \in I$ with $t \leq \tau < t + \delta(\varepsilon, t)$ it holds that $|f(\tau) - f(t)| < \varepsilon$. A continuous function $f : [0,1] \to \Re$ is called piecewise $C^1$ on $[0,1]$ and we write $f \in PC^1([0,1])$, if the following properties hold: (i) for every $x \in [0,1)$ the limit $\lim_{h \to 0^+} \left( h^{-1} (f(x+h) - f(x)) \right)$ exists and is finite, (ii) for every $x \in (0,1]$ the limit $\lim_{h \to 0^-} \left( h^{-1} (f(x+h) - f(x)) \right)$ exists and is finite, (iii) there exists a finite set $I \subset (0,1)$ where $f'(x) = \lim_{h \to 0^-} \left( h^{-1} (f(x+h) - f(x)) \right) = \lim_{h \to 0^+} \left( h^{-1} (f(x+h) - f(x)) \right)$ holds for $x \in (0,1) \setminus I$, and (iv) the mapping $(0,1) \setminus I \ni x \to f'(x) \in \Re$ is continuous.

* Let $r \in C^0([0,1]; (0,+\infty))$ be given. $L^2_r(0,1)$ denotes the equivalence class of measurable functions $f : [0,1] \to \Re$ for which $\|f\|_r = \left( \int_0^1 r(x) |f(x)|^2 dx \right)^{1/2} < +\infty$. $L^2_r(0,1)$ is a Hilbert space with inner product $\langle f, g \rangle = \int_0^1 r(x) f(x) g(x) dx$. When $r(x) \equiv 1$, we use the notation $L^2(0,1)$ for the standard space of square-integrable functions and $\|f\| = \left( \int_0^1 |f(x)|^2 dx \right)^{1/2} < +\infty$ for $f \in L^2(0,1)$.

* Let $u : \Re_+ \times [0,1] \to \Re$ be given. We use the notation $u[t]$ to denote the profile at certain $t \geq 0$, i.e., $(u[t])(x) = u(t, x)$ for all $x \in [0,1]$. When $u(t, x)$ is differentiable with respect to $x \in [0,1]$, we use the notation $u'(t, x)$ for the derivative of $u$ with respect to $x \in [0,1]$, i.e., $u'(t, x) = \frac{\partial u}{\partial x}(t, x)$. For an interval $I \subseteq \Re_+$, the space $C^0(I; L^2_r(0,1))$ is the space of continuous mappings $I \ni t \to u[t] \in L^2_r(0,1)$.

* $H^1(0,1)$ denotes the Sobolev space of continuous functions on $[0,1]$ with measurable, square integrable derivative. $H^2(0,1)$ denotes the Sobolev space of continuously differentiable functions on $[0,1]$ with measurable, square integrable second derivative.

## 2. Main Results

Consider the Sturm-Liouville (SL) operator $A : D \to L^2_r(0,1)$ defined by

$$(Af)(x) = -\frac{1}{r(x)} \frac{d}{dx}\left( p(x) \frac{df}{dx}(x) \right) + \frac{q(x)}{r(x)} f(x), \text{ for all } f \in D \text{ and } x \in (0,1) \tag{2.1}$$

where $p \in C^1([0,1]; (0,+\infty))$, $r \in C^0([0,1]; (0,+\infty))$, $q \in C^0([0,1]; \Re)$ and $D \subseteq H^2(0,1)$ is the set of all functions $f : [0,1] \to \Re$ for which

$$b_1 f(0) + b_2 f'(0) = a_1 f(1) + a_2 f'(1) = 0 \tag{2.2}$$

where $a_1, a_2, b_1, b_2$ are real constants with $|a_1| + |a_2| > 0$, $|b_1| + |b_2| > 0$. It is well-known (Chapter 11 in [3] and pages 498-505 in [34]) that all eigenvalues of the SL operator $A : D \to L^2_r(0,1)$, defined by (2.1), (2.2) are real. The eigenvalues form an infinite, increasing sequence $\lambda_1 < \lambda_2 < \ldots < \lambda_n < \ldots$ with $\lim_{n \to \infty} (\lambda_n) = +\infty$ and to each eigenvalue $\lambda_n \in \Re$ ($n = 1, 2, \ldots$) corresponds exactly one eigenfunction $\phi_n \in C^2([0,1]; \Re)$ that satisfies $A\phi_n = \lambda_n \phi_n$ and $b_1 \phi_n(0) + b_2 \phi'_n(0) = a_1 \phi_n(1) + a_2 \phi'_n(1) = 0$. Moreover, the eigenfunctions form an orthonormal basis of $L^2_r(0,1)$.

In the present work, we use the following assumption for the SL operator $A: D \to L_r^2(0,1)$ defined by (2.1), (2.2), where $a_1, a_2, b_1, b_2$ are real constants with $|a_1|+|a_2|>0$, $|b_1|+|b_2|>0$.

**(H):** *The SL operator $A: D \to L_r^2(0,1)$ defined by (2.1), (2.2), where $a_1, a_2, b_1, b_2$ are real constants with $|a_1|+|a_2|>0$, $|b_1|+|b_2|>0$, satisfies*

$$\sum_{n=N}^{\infty} \lambda_n^{-1} \max_{0 \leq x \leq 1}(|\phi_n(x)|) < +\infty, \text{ for certain } N>0 \text{ with } \lambda_N > 0 \qquad (2.3)$$

It is important to notice that the validity of Assumption (H) can be verified without knowledge of eigenvalues and the eigenfunctions of the SL operator $A$ (see [35]).

We next consider the following system

$$\frac{\partial u}{\partial t}(t,x) - \frac{1}{r(x)} \frac{\partial}{\partial x}\left(p(x)\frac{\partial u}{\partial x}(t,x)\right) + \frac{q(x)}{r(x)} u(t,x) = f(t,x), \quad x \in (0,1) \qquad (2.4)$$

$$b_1 u(t,0) + b_2 \frac{\partial u}{\partial x}(t,0) - d_0(t) = a_1 u(t,1) + a_2 \frac{\partial u}{\partial x}(t,1) - d_1(t) = 0, \qquad (2.5)$$

where $u[t]$ is the state and $f(t,x), d_0(t), d_1(t)$ are disturbance inputs. The following definition gives classes of functions in which disturbance inputs may belong.

**Definition 2.1:** *A function $f: \Re_+ \times [0,1] \to \Re$ is of class GD, if there exists an increasing sequence of times $\{\tau_i \geq 0, i = 0,1,2,...\}$ with $\tau_0 = 0$, $\lim_{i \to +\infty}(\tau_i) = +\infty$ with the following property: for every $i \geq 1$ there exist functions $f_{1,i}, f_{2,i} \in C^0((\tau_{i-1}, \tau_i) \times [0,1])$, $a_i \in PC^1([0,1])$, $b_i \in C^0([0,1])$, $c_i \in D \cap C^2([0,1])$, with $f_{1,i}[t] \in PC^1([0,1])$, $f_{2,i}[t] \in D \cap C^2([0,1])$, $f[t] = f_1[t] + f_2[t]$ for all*

$$t \in (\tau_{i-1}, \tau_i), \sup_{t \in (\tau_{i-1}, \tau_i)}\left(\int_0^1 r(x)\left(|f_{1,i}(t,x)|^2 + |f_{2,i}(t,x)|^2 + |(Af_{2,i}[t])(x)|^2\right)dx\right) < +\infty, (\tau_{i-1}, \tau_i) \times [0,1] \ni (t,x) \to \frac{\partial f_{1,i}}{\partial t}(t,x)$$

*is continuous, $\lim_{t \to \tau_i^-}(f_{1,i}(t,x)) = a_i(x)$, $\lim_{t \to \tau_i^-}\left(\frac{\partial f_{1,i}}{\partial t}(t,x)\right) = b_i(x)$, $\lim_{t \to \tau_i^-}(f_{2,i}(t,x)) = c_i(x)$ for all $x \in [0,1]$.*

*A right continuous function $d: \Re_+ \to \Re$ is of class GB, if there exists an increasing sequence of times $\{\tau_i \geq 0, i = 0,1,2,...\}$ with $\tau_0 = 0$, $\lim_{i \to +\infty}(\tau_i) = +\infty$ such that:*

(i) $d \in C^2(I)$, where $I = \Re_+ \setminus \{\tau_i \geq 0, i = 0,1,2,...\}$,
(ii) *for every $\tau_i > 0$ the left limits of $d(t), \dot{d}(t), \ddot{d}(t)$ when $t$ tends to $\tau_i$ are finite,*
(iii) $\sup_{t \in (\tau_i, \tau_{i+1})}(|\dot{d}(t)|) < +\infty$ *for $i = 0,1,2,...$ are finite.*

The first result clarifies the properties of the solution of (2.4) and (2.5) under the presence of the disturbance inputs $f(t,x), d_0(t), d_1(t)$.

**Theorem 2.2:** *Consider the SL operator $A: D \to L_r^2(0,1)$ defined by (2.1), (2.2), where $a_1, a_2, b_1, b_2$ are real constants with $|a_1|+|a_2|>0$, $|b_1|+|b_2|>0$, under Assumption (H). Let $f \in GD$, $d_0, d_1 \in GB$ be given functions and let $\{\tau_i \geq 0, i = 0,1,2,...\}$ be the increasing sequence of times with $\tau_0 = 0$, $\lim_{i \to +\infty}(\tau_i) = +\infty$ involved in Definition 2.1. Then for every $u_0 \in L_r^2(0,1)$ there exists a unique mapping $u \in C^0(\Re_+; L_r^2(0,1))$ with $u \in C^1(I \times [0,1])$ satisfying $u[t] \in C^2([0,1])$ for all $t > 0$, $\lim_{t \to \tau_i^-}(u(t,x)) = u(\tau_i, x)$, $\lim_{t \to \tau_i^-}\left(\frac{\partial u}{\partial t}(t,x)\right) = -(Au[\tau_i])(x) + \lim_{t \to \tau_i^-}(f(t,x))$, $\lim_{t \to \tau_i^-}\left(\frac{\partial u}{\partial z}(t,x)\right) = \frac{\partial u}{\partial z}(\tau_i, x)$ for all $i \geq 1$, $u(0,x) = u_0(x)$ for all $x \in [0,1]$, and*

$$\frac{\partial u}{\partial t}(t,x) - \frac{1}{r(x)}\frac{\partial}{\partial x}\left(p(x)\frac{\partial u}{\partial x}(t,x)\right) + \frac{q(x)}{r(x)}u(t,x) = f(t,x), \text{ for all } (t,x) \in I \times (0,1) \quad (2.6)$$

$$b_1 u(t,0) + b_2 \frac{\partial u}{\partial x}(t,0) - d_0(t) = a_1 u(t,1) + a_2 \frac{\partial u}{\partial x}(t,1) - d_1(t) = 0, \text{ for all } t \in I \quad (2.7)$$

where $I = \Re_+ \setminus \{\tau_i \geq 0, i = 0,1,2,...\}$.

In what follows, for any given $u_0 \in L^2_r(0,1)$, $f \in GD$, $d_0, d_1 \in GB$, the mapping $u \in C^0(\Re_+; L^2_r(0,1))$ with $u \in C^1(I \times [0,1])$ satisfying $u[t] \in C^2([0,1])$ for all $t > 0$, $\lim_{t \to \tau_i^-}(u(t,x)) = u(\tau_i, x)$, $\lim_{t \to \tau_i^-}\left(\frac{\partial u}{\partial t}(t,x)\right) = -(Au[\tau_i])(x) + \lim_{t \to \tau_i^-}(f(t,x))$, $\lim_{t \to \tau_i^-}\left(\frac{\partial u}{\partial x}(t,x)\right) = \frac{\partial u}{\partial x}(\tau_i, x)$ for all $i \geq 1$, $u(0,x) = u_0(x)$ for all $x \in [0,1]$, and (2.6), (2.7), is called the (unique) *solution of the evolution equation (2.4) with (2.5) and initial condition* $u_0 \in L^2_r(0,1)$ *corresponding to inputs* $f \in GD$, $d_0, d_1 \in GB$.

**Definition 2.3:** *Consider the SL operator* $A: D \to L^2_r(0,1)$ *defined by (2.1), (2.2), where* $a_1, a_2, b_1, b_2$ *are real constants with* $|a_1| + |a_2| > 0$, $|b_1| + |b_2| > 0$, *under Assumption (H). The operator* $A: D \to L^2_r(0,1)$ *is called Exponentially Stable (ES) if* $\lambda_1 > 0$.

Our second main result provides decay estimates of the solution of (2.4), (2.5) in the norm of $L^2_r(0,1)$.

**Theorem 2.4:** *Suppose that the SL operator* $A: D \to L^2_r(0,1)$ *defined by (2.1), (2.2) is ES. Then for every* $u_0 \in L^2_r(0,1)$, $f \in GD$, $d_0, d_1 \in GB$, *the unique solution* $u \in C^0(\Re_+; L^2_r(0,1))$ *of the evolution equation (2.4) with (2.5) and initial condition* $u_0 \in L^2_r(0,1)$ *satisfies the following estimate for all* $\sigma \in [0, \lambda_1)$ *and* $t > 0$:

$$\|u[t]\|_r \leq \exp(-\lambda_1 t)\|u_0\|_r + \frac{\lambda_1}{\lambda_1 - \sigma} C_0 \sup_{0 < s < t}\left(|d_0(s)|\exp(-\sigma(t-s))\right)$$
$$+ \frac{\lambda_1}{\lambda_1 - \sigma} C_1 \sup_{0 < s < t}\left(|d_1(s)|\exp(-\sigma(t-s))\right) + \frac{1}{\lambda_1 - \sigma} \sup_{0 < s < t}\left(\|f[s]\|_r \exp(-\sigma(t-s))\right) \quad (2.8)$$

where

$$C_0 := \frac{p(0)}{b_1^2 + b_2^2}\sqrt{\sum_{n=1}^{\infty}\frac{1}{\lambda_n^2}\left|b_1 \frac{d\phi_n}{dz}(0) - b_2\phi_n(0)\right|^2} = \frac{1}{\sqrt{b_1^2 + b_2^2}}\|\tilde{u}\|_r, \quad (2.9)$$

$$C_1 := \frac{p(1)}{a_1^2 + a_2^2}\sqrt{\sum_{n=1}^{\infty}\frac{1}{\lambda_n^2}\left|a_2\phi_n(1) - a_1\frac{d\phi_n}{dz}(1)\right|^2} = \frac{1}{\sqrt{a_1^2 + a_2^2}}\|\bar{u}\|_r, \quad (2.10)$$

$\tilde{u} \in C^2([0,1])$ *is the unique solution of the boundary value problem* $(p(x)\tilde{u}'(x))' - q(x)\tilde{u}(x) = 0$ *for* $x \in [0,1]$ *with* $b_1\tilde{u}(0) + b_2\tilde{u}'(0) = \sqrt{b_1^2 + b_2^2}$, $a_1\tilde{u}(1) + a_2\tilde{u}'(1) = 0$ *and* $\bar{u} \in C^2([0,1])$ *is the unique solution of the boundary value problem* $(p(x)\bar{u}'(x))' - q(x)\bar{u}(x) = 0$ *for* $x \in [0,1]$ *with* $b_1\bar{u}(0) + b_2\bar{u}'(0) = 0$ *and* $a_1\bar{u}(1) + a_2\bar{u}'(1) = \sqrt{a_1^2 + a_2^2}$.

**Remark 2.5:** Theorem 2.4 generalizes in many ways the result of Theorem 2.3 in [23]. Since Theorem 2.3 in [23] considers inputs for which the solution $u$ is continuous on $\Re_+ \times [0,1]$ and $C^1$ on $(0, +\infty) \times [0,1]$, it follows that no discontinuous inputs are allowed. On the other hand, Theorem 2.4 allows discontinuous inputs. Another difference between Theorem 2.3 in [23] and Theorem 2.4 is the obtained decay estimate. While Theorem 2.3 in [23] provides an decay estimate which involves only terms of the form $\sup_{0\leq s\leq t}(|d_0(s)|\exp(-\sigma(t-s)))$, $\sup_{0\leq s\leq t}(|d_1(s)|\exp(-\sigma(t-s)))$ with $\sigma=0$, Theorem 2.4 allows positive values for $\sigma$. The difference is important, because Theorem 2.4 provides a "fading memory estimate" (see [19]), which can be directly used for small-gain analysis. The term describing the effect of the initial condition in (2.8), namely the term $\exp(-\lambda_1 t)\|u_0\|_r$, is less conservative than the one used in Theorem 2.3 in [23]. Finally, Theorem 2.3 in [23] used the infinity norm of the distributed input $f$, while Theorem 2.4 uses the weighted 2-norm of the distributed input $f$. The proofs of Theorem 2.3 in [23] and Theorem 2.4 are very different (due to the fact that Theorem 2.4 does not impose the demanding regularity requirements for the solution of (2.4) and (2.5) in [23]).

It is clear that formulas (2.9), (2.10) allow us to calculate the gain coefficients in (2.8) without requiring knowledge of the eigenvalues and the eigenfunctions of the SL operator $A:D\to L_r^2(0,1)$. However, we still need to know a positive lower bound of the principal eigenvalue $\lambda_1$. The following proposition provides the means to avoid the calculation of $\lambda_1$, as well as the means to verify that the SL operator $A:D\to L_r^2(0,1)$ defined by (2.1), (2.2) is ES. It deals with the special case $p(x)\equiv p$, $r(x)\equiv 1$, because, when $p\in C^2([0,1];(0,+\infty))$, $r\in C^2([0,1];(0,+\infty))$, using a so-called gauge transformation, it is possible to convert system (2.4), (2.5) into one with constant diffusion and zero advection terms. More specifically, this is achieved by the coordinate change:

$$\xi=\sqrt{\varepsilon}\int_0^x \sqrt{\frac{p(s)}{r(s)}}ds\ ,\ \text{where}\ \varepsilon=\left(\int_0^1\sqrt{\frac{r(s)}{p(s)}}ds\right)^{-2}\ \text{and}\ U(t,\xi)=(r(x)p(x))^{1/4}u(t,x).$$

**Proposition 2.6:** *Consider the SL operator $A:D\to L^2(0,1)$ defined by (2.1), (2.2), where $a_1,a_2,b_1,b_2$ are real constants with $|a_1|+|a_2|>0$, $|b_1|+|b_2|>0$, with $p(x)\equiv p$, $r(x)\equiv 1$, under Assumption (H). Suppose that there exist constants $\varepsilon_1,\varepsilon_2>0$, $\lambda\in[0,1]$ and a function $g\in C^2([0,1];(0,+\infty))$ such that*

$$g'(1)-2q_1g(1)-2R^{-1}(1-\lambda)(1+\varepsilon_1)\geq 0 \tag{2.11}$$
$$2q_0g(0)-g'(0)-2R^{-1}\lambda(1+\varepsilon_0)\geq 0$$

*and*

$$2R^{-1}+2p^{-1}g(x)q(x)>g''(x),\ \text{for all}\ x\in[0,1] \tag{2.12}$$

*where*

$$R:=\int_0^1\left(\lambda(1+\varepsilon_0^{-1})\int_0^z\frac{ds}{g(s)}+(1-\lambda)(1+\varepsilon_1^{-1})\int_z^1\frac{ds}{g(s)}\right)dz \tag{2.13}$$

$$q_0=+\infty\ \text{if}\ b_2=0\ \text{and}\ q_0=-b_1/b_2\ \text{if}\ b_2\neq 0 \tag{2.14}$$
$$q_1=-\infty\ \text{if}\ a_2=0\ \text{and}\ q_1=-a_1/a_2\ \text{if}\ a_2\neq 0 \tag{2.15}$$

*Then the operator $A:D\to L^2(0,1)$ is ES. Moreover, $\lambda_1\geq\dfrac{\min\{2g(x)q(x)+2pR^{-1}-pg''(x):x\in[0,1]\}}{2\max\{g(x):x\in[0,1]\}}>0$.*

In order to obtain decay estimates in the $H^1$ norm we need to focus our attention to more specific cases. The following result provides decay estimates for the Dirichlet case in the $L^2$ norm of the solution as well as in the $L^2$ norm of the spatial derivative of the solution.

**Theorem 2.7 (Dirichlet BCs-no boundary disturbances):** *Consider the SL operator $A:D\to L^2(0,1)$, defined by (2.1), (2.2) with $p(x)\equiv p$, $r(x)\equiv 1$, $q(x)\equiv q$, $a_2=b_2=0$, $a_1=b_1=1$, where $q>-p\pi^2$. Then for every $u_0\in H^1(0,1)$ with $u_0(0)=u_0(1)=0$ and $f\in GD$ with $f[t]\in C^1([0,1])$ for all $t\geq 0$, $f'\in GD$ and $f(t,1)$, $f(t,0)$ being of class GB, the unique solution $u\in C^0(\mathfrak{R}_+;L^2(0,1))$ of the* *evolution equation (2.4) with (2.5), $d_0(t) = d_1(1) \equiv 0$ and initial condition $u_0 \in H^1(0,1)$ satisfies the following estimates for all $\sigma \in [0, p\pi^2 + q)$, $t > 0$:*

$$\|u[t]\| \leq \exp\left(-\left(p\pi^2 + q\right)t\right)\|u_0\| + \frac{1}{p\pi^2 + q - \sigma} \sup_{0 < s < t}\left(\|f[s]\|\exp(-\sigma(t-s))\right) \tag{2.16}$$

$$\|u'[t]\| \leq \exp\left(-\left(p\pi^2 + q\right)t\right)\|u'_0\| + \left(\frac{p\pi^2 + q}{p\pi^2 + q - \sigma}\right)p^{-1}h(q/p) \sup_{0 < s < t}\left(\|f[s]\|\exp(-\sigma(t-s))\right) \tag{2.17}$$

*where*

$$h(s) := \begin{cases} \dfrac{\sqrt{\exp(2\sqrt{s}) - \exp(-2\sqrt{s}) - 4\sqrt{s}}}{2(s)^{1/4}\left(\exp(\sqrt{s}) - \exp(-\sqrt{s})\right)} & \text{if } s > 0 \\ \dfrac{\sqrt{6}}{6} & \text{if } s = 0 \\ \dfrac{\sqrt{4\sqrt{-s} - 2\sin(2\sqrt{-s})}}{4(-s)^{1/4}\sin(\sqrt{-s})} & \text{if } -\pi^2 < s < 0 \end{cases} \tag{2.18}$$

In order to obtain decay estimates in the $H^1$ norm for the case with Dirichlet BC at 0 and Robin (or Neumann) BC at 1 we need the following definition.

**Definition 2.8:** *Consider the SL operator $A: D \to L^2_r(0,1)$ defined by (2.1), (2.2), where $a_1, a_2, b_1, b_2$ are real constants with $q \in C^2([0,1];(0,+\infty))$, $b_2 = 0$, $b_1 = a_2 = 1$, $p(x) \equiv p$, $r(x) \equiv 1$, under Assumption (H). The "derived operator" of $A$ is the SL operator $A': D' \to L^2(0,1)$ defined by*

$$(A'f)(x) = -pf''(x) + (q(x) + 2a_1 p)f(x), \text{ for all } f \in D' \text{ and } x \in (0,1) \tag{2.19}$$

*where $D' \subseteq H^2(0,1)$ is the set of all functions $f : [0,1] \to \Re$ for which*

$$f'(0) = f(1) = 0 \tag{2.20}$$

**Theorem 2.9 (Dirichlet BC at 0 with no boundary disturbance-Robin BC at 1 with boundary disturbance):** *Suppose that the SL operator $A: D \to L^2(0,1)$ with $q \in C^2([0,1];(0,+\infty))$, $p(x) \equiv p$, $r(x) \equiv 1$, $b_2 = 0$, $b_1 = a_2 = 1$ and its derived operator $A': D' \to L^2(0,1)$ defined by (2.19), (2.20) are ES. Then there exist constants $M, \Theta_1, \Theta_2, \Theta_3, \gamma_0, \gamma_1, \gamma_2, \sigma > 0$ such that for every $u_0 \in H^1(0,1)$ with $u_0(0) = 0$, $d_1 \in GB$, $f \in GD$ with $f[t] \in C^1([0,1])$ for all $t \geq 0$, $f' \in GD$ and $f(t,0)$ being of class GB, the unique solution $u \in C^0(\Re_+; L^2(0,1))$ of the evolution equation (2.4) with (2.5), $d_0(t) \equiv 0$ and initial condition $u_0 \in H^1(0,1)$ satisfies the following estimates for all $t > 0$:*

$$\|u'[t]\| \leq \left(\|u'[0]\| + M\|u[0]\|\right)\exp(-\sigma t) + \Theta_1 \sup_{0 \leq s \leq t}\left(\|f[s]\|\exp(-\sigma(t-s))\right) + \Theta_2 \sup_{0 < s < t}\left(\|f'[s]\|\exp(-\sigma(t-s))\right)$$
$$+ \gamma_0 \sup_{0 < s < t}\left(|f(s,0)|\exp(-\sigma(t-s))\right) + \gamma_1 \sup_{0 < s < t}\left(|d_1(s)|\exp(-\sigma(t-s))\right) \tag{2.21}$$

$$\|u[t]\| \leq \exp(-\sigma t)\|u_0\| + \Theta_3 \sup_{0 < s < t}\left(\|f[s]\|\exp(-\sigma(t-s))\right) + \gamma_2 \sup_{0 < s < t}\left(|d_1(s)|\exp(-\sigma(t-s))\right) \tag{2.22}$$

**Remark 2.10:** The proof of Theorem 2.9 provides explicit estimates for all constants $M, \Theta_1, \Theta_2, \gamma_0, \gamma_1, \sigma > 0$ appearing in the right hand side of (2.21). More specifically, we get

$$M = 2|a_1| + \frac{\max_{0 \leq x \leq 1}\left(|q'(x) - 2a_1^2 px|\right)}{\mu_1 - \sigma}$$

$$\gamma_0 = \frac{\mu_1}{\mu_1 - \sigma}\tilde{C}_0 p^{-1} \quad , \quad \gamma_1 = \frac{\mu_1}{\mu_1 - \sigma}\tilde{C}_1 + \frac{\lambda_1}{\lambda_1 - \sigma}C_1\left(\frac{\max_{0 \leq x \leq 1}\left(|q'(x) - 2a_1^2 px|\right)}{\mu_1 - \sigma} + |a_1|\right) \quad , \quad \gamma_2 = \frac{\lambda_1}{\lambda_1 - \sigma}C_1$$

$$\Theta_1 = \frac{|a_1|}{\mu_1 - \sigma} + \frac{|a_1|}{\lambda_1 - \sigma} + \lambda_1 \frac{\max_{0 \leq x \leq 1}\left(|q'(x) - 2a_1^2 px|\right)}{(\mu_1 - \sigma)(\lambda_1 - \sigma)} \quad , \quad \Theta_2 = \frac{1}{\mu_1 - \sigma} \quad , \quad \Theta_3 = \frac{1}{\lambda_1 - \sigma}$$

where $\sigma \in [0, \min(\mu_1, \lambda_1))$ is arbitrary, $C_0 := \|\tilde{u}\|$, $C_1 := \frac{1}{\sqrt{a_1^2 + 1}}\|\tilde{u}\|$, $\tilde{u} \in C^2([0,1])$ is the unique solution of the boundary value problem $p\tilde{u}''(z) - q(x)\tilde{u}(x) = 0$ for $x \in [0,1]$ with $\tilde{u}(0) = 1$, $a_1\tilde{u}(1) + \tilde{u}'(1) = 0$, $\bar{u} \in C^2([0,1])$ is the unique solution of the boundary value problem $p\bar{u}''(x) - q(x)\bar{u}(z) = 0$ for $x \in [0,1]$ with $\bar{u}(0) = 0$ and $a_1\bar{u}(1) + \bar{u}'(1) = \sqrt{a_1^2 + 1}$, $0 < \mu_1 < \mu_2 < \ldots < \mu_n < \ldots$ with $\lim_{n \to \infty}(\mu_n) = +\infty$ are the eigenvalues of the SL operator $A'$, $\tilde{C}_0 := \|\tilde{v}\|$, $\tilde{C}_1 := \|\tilde{v}\|$, $\tilde{v} \in C^2([0,1])$ is the solution of the boundary value problem $p\tilde{v}''(x) - (q(x) + 2a_1 p)\tilde{v}(x) = 0$ for $x \in [0,1]$ with $\tilde{v}'(0) = 1$, $\tilde{v}(1) = 0$ and $\bar{v} \in C^2([0,1])$ is the solution of the boundary value problem $p\bar{v}''(x) - (q(x) + 2a_1 p)\bar{v}(x) = 0$ for $x \in [0,1]$ with $\bar{v}'(0) = 0$ and $\bar{v}(1) = 1$.

The proofs of Theorem 2.8 and Theorem 2.9 rely on the following technical results, which are of independent interest. Both results show that the spatial derivative of the solution of the evolution equation (2.4) with (2.5) is determined by solving a specific evolution equation with specific boundary disturbances (even if boundary disturbances were absent in the original evolution equation (2.4), (2.5)).

**Proposition 2.11:** *Consider the SL operator $A : D \to L^2(0,1)$ defined by (2.1), (2.2) with $p(x) \equiv p$, $r(x) \equiv 1$, $q(x) \equiv q$, $a_2 = b_2 = 0$, $a_1 = b_1 = 1$. Let $u_0 \in H^1(0,1)$ with $u_0(0) = u_0(1) = 0$, $f \in GD$ with $f[t] \in C^1([0,1])$ for all $t \geq 0$, $f' \in GD$ and $f(t,1)$, $f(t,0)$ being of class $GB$, be given functions. Consider the solution $u \in C^0(\Re_+; L^2(0,1))$ of the evolution equation (2.4) with (2.5), $d_0(t) = d_1(1) \equiv 0$ and initial condition $u_0 \in H^1(0,1)$ with $u_0(0) = u_0(1) = 0$, corresponding to input $f \in GD$. Consider also the solution $v \in C^0(\Re_+; L^2(0,1))$ of the initial-boundary value problem*

$$\frac{\partial v}{\partial t}(t,x) = p\frac{\partial^2 v}{\partial x^2}(t,x) - qv(t,x) + f'(t,x) \tag{2.23}$$

$$\frac{\partial v}{\partial x}(t,0) + p^{-1}f(t,0) = \frac{\partial v}{\partial x}(t,1) + p^{-1}f(t,1) = 0 \tag{2.24}$$

*with initial condition $v_0 = u_0'$. Then the following equations hold for all $(t,x) \in \Re_+ \times [0,1]$:*

$$v(t,x) = \frac{\partial u}{\partial x}(t,x) \tag{2.25}$$

$$\int_0^1 v(t,x)dx = 0 \tag{2.26}$$

**Proposition 2.12:** *Consider the SL operator $A : D \to L^2(0,1)$ defined by (2.1), (2.2) with $p(x) \equiv p$, $r(x) \equiv 1$, $b_2 = 0$, $a_2 = b_1 = 1$. Suppose that the SL operator $A' : D' \to L^2(0,1)$ defined by (2.19), (2.20), satisfies assumption (H). Let $u_0 \in H^1(0,1)$ with $u_0(0) = 0$, $d_1 \in GB$, $f \in GD$ with $f[t] \in C^1([0,1])$ for all $t \geq 0$, $f' \in GD$ and $f(t,0)$ being of class $GB$, be given functions. Consider the solution $u \in C^0(\Re_+; L^2(0,1))$ of the evolution equation (2.4) with (2.5), $d_0(t) \equiv 0$ and initial condition* $u_0 \in H^1(0,1)$, *corresponding to inputs* $f \in GD$, $d_1 \in GB$. *Consider also the solution* $v \in C^0(\Re_+; L^2(0,1))$ *of the initial-boundary value problem*

$$\frac{\partial v}{\partial t}(t,x) = p\frac{\partial^2 v}{\partial x^2}(t,x) - (q(x) + 2a_1 p)v(t,x) - (q'(x) - 2a_1^2 px)u(t,z) + f'(t,x) + a_1 x f(t,x) \quad (2.27)$$

$$\frac{\partial v}{\partial x}(t,0) + p^{-1}f(t,0) = v(t,1) - d_1(t) = 0 \quad (2.28)$$

*with initial condition* $v_0(x) = u_0'(x) + a_1 x u_0(x)$ *for* $x \in [0,1]$. *Then the following equation holds for all* $(t,x) \in \Re_+ \times [0,1]$:

$$v(t,x) = \frac{\partial u}{\partial x}(t,x) + a_1 x u(t,x) \quad (2.29)$$

The following example shows how easily the obtained results can be applied to the study of heat transfer phenomena.

**Example 2.13:** Consider a solid bar of length $L > 0$ and its temperature $T(t,x)$ at time $t \geq 0$ and position $z \in [0, L]$. The temperature of the bar is kept constant at $z = 0$, i.e..,

$$T(t,0) = T_0, \text{ for } t \geq 0 \quad (2.30)$$

while at $z = 1$ the bar is in contact with air. The air temperature $T_{air}(t)$ is subject to variation around a nominal temperature $T_{nom}$, i.e.,

$$T_{air}(t) = T_{nom} + d(t), \text{ for } t \geq 0 \quad (2.31)$$

Applying Newton's law of cooling and Fourier's law of heat conduction, we get

$$-k\frac{\partial T}{\partial z}(t,L) = h(T(t,L) - T_{air}(t)) \text{ for } t \geq 0 \quad (2.32)$$

where $h > 0$ is the heat transfer coefficient of air and $k > 0$ is the thermal conductivity of the solid. Taking into account (2.30), (2.31), (2.32) and using the dimensionless position $x = z/L$, we obtain the following evolution problem:

$$\frac{\partial T}{\partial t}(t,x) = p\frac{\partial^2 T}{\partial x^2}(t,x), \text{ for } x \in (0,1) \quad (2.33)$$

$$T(t,0) - T_0 = \frac{\partial T}{\partial x}(t,1) + a_1 T(t,1) - a_1 T_{nom} - a_1 d(t) = 0, \quad (2.34)$$

where $p, a_1 > 0$ are constants. Using the dimensionless deviation variable $u(t,x)$ around the equilibrium temperature profile $T(x) = T_0 + \frac{a_1}{1+a_1}(T_{nom} - T_0)x$, i.e., defining

$$u(t,x) = \frac{T(t,x)}{T_0} - 1 - \frac{a_1}{1+a_1}\left(\frac{T_{nom}}{T_0} - 1\right)x \quad (2.35)$$

we get the following evolution problem:

$$\frac{\partial u}{\partial t}(t,x) = p\frac{\partial^2 u}{\partial x^2}(t,x), \text{ for } x \in (0,1) \quad (2.36)$$

$$u(t,0) = \frac{\partial u}{\partial x}(t,1) + a_1 u(t,1) - d_1(t) = 0, \quad (2.37)$$

where $d_1(t) = a_1 \dfrac{d(t)}{T_0}$.

The goal is the estimation of the effect of the disturbance $d_1(t)$ to the temperature profile of the bar. The evolution of $u$ may be studied using the results in [21,22,23]. In this case, we obtain estimates for initial conditions $u_0 \in C^2([0,1])$ and disturbances $d_1 \in C^2(\Re_+)$ with $u_0(0) = u_0'(1) + a_1 u_0(1) - d_1(0) = 0$.

Let $\theta \in \left(\dfrac{\pi}{2}, \pi\right)$ be the unique solution of the equation $\tan(\theta) = -a_1^{-1}\theta$. It follows from Theorem 2.3 in [23] or Theorem 2.2 in [22] (by performing all relevant computations) that the following estimate holds for all $t \geq 0$:

$$\|u[t]\| \leq \sqrt{\dfrac{\exp(-p\theta^2 t)}{2 - \exp(-p\theta^2 t)}} \|u_0\| + \dfrac{\sqrt{3}}{3(1+a_1)} \max_{0 \leq s \leq t}(|d_1(s)|) \qquad (2.38)$$

Pick $\theta_0 \in \left(\dfrac{\pi}{2}, \theta\right)$, $\varphi \in (0, \theta - \theta_0)$ and notice that Assumption (H4) in [23] holds with $\eta(x) = \sin(\theta_0 x + \varphi)$. It follows from Theorem 2.2 in [23] that the following estimate holds for all $t \geq 0$, $\theta_0 \in \left(\dfrac{\pi}{2}, \theta\right)$ and $\varphi \in (0, \theta - \theta_0)$:

$$\max_{0 \leq x \leq 1}\left(\dfrac{|u(t,x)|}{\sin(\theta_0 x + \varphi)}\right) \leq \max\left(\exp(-p\theta_0^2 t)\max_{0 \leq x \leq 1}\left(\dfrac{|u_0(x)|}{\sin(\theta_0 x + \varphi)}\right), \dfrac{\max_{0 \leq s \leq t}(|d_1(s)|)}{a_1 \sin(\theta_0 + \varphi) + \theta_0 \cos(\theta_0 + \varphi)}\right) \qquad (2.39)$$

Applying the results of the present work, we get different results. Theorem 2.4 can be applied to initial conditions $u_0 \in L^2(0,1)$ and disturbances $d_1 \in GB$ and provides the following estimate which holds for all $t > 0$ and $\sigma \in [0, p\theta^2)$:

$$\|u[t]\| \leq \exp(-p\theta^2 t)\|u_0\| + \dfrac{p\theta^2}{p\theta^2 - \sigma} \dfrac{\sqrt{3}}{3(1+a_1)} \sup_{0 < s < t}(|d_1(s)|\exp(-\sigma(t-s))) \qquad (2.40)$$

This estimate should be compared with estimate (2.38). Moreover, by performing all relevant computations, we can verify that Theorem 2.9 can be also applied. Using Theorem 2.9 and Remark 2.10, we are in a position to conclude that for initial conditions $u_0 \in H^1(0,1)$ with $u_0(0) = 0$ and disturbances $d_1 \in GB$, the following estimate which holds for all $t > 0$ and $\sigma \in \left[0, \min\left(p\dfrac{\pi^2 + 8a_1}{4}, p\theta^2\right)\right]$:

$$\|u'[t]\| \leq (\|u'[0]\| + M\|u[0]\|)\exp(-\sigma t) + \gamma_1 \sup_{0 < s < t}(|d_1(s)|\exp(-\sigma(t-s))) \qquad (2.41)$$

where

$$M := 2a_1 + \dfrac{8a_1^2 p}{p\pi^2 + 8pa_1 - 4\sigma}$$

$$\gamma_1 := \dfrac{p\pi^2 + 8pa_1}{p\pi^2 + 8pa_1 - 4\sigma}\tilde{C}_1 + \dfrac{p\theta^2}{p\theta^2 - \sigma} \dfrac{\sqrt{3}a_1}{3(1+a_1)} \dfrac{p\pi^2 + 16pa_1 - 4\sigma}{p\pi^2 + 8pa_1 - 4\sigma}$$

$$\tilde{C}_1 := \dfrac{1}{\exp(\sqrt{2a_1}) + \exp(-\sqrt{2a_1})}\sqrt{2 + \dfrac{\exp(2\sqrt{2a_1}) - \exp(-2\sqrt{2a_1})}{2\sqrt{2a_1}}}$$

Estimates (2.39), (2.40) and (2.41) may be used in a straightforward way in order to obtain quantitative results for the temperature of the bar. ◁

# 3. Applications to 1-D Nonlocal Parabolic PDEs

This section provides two examples of 1-D parabolic PDEs with nonlocal terms. The nonlocal terms may appear either in the BCs or in the PDE. We are not dealing with existence/uniqueness issues for the PDEs; instead we show how the main results of the present work can be used directly in order to derive exponential stability estimates for the solution.

**Example 3.1:** Let $p>0$, $q\in\Re$ with $q>-p\pi^2/4$ and consider the system

$$\frac{\partial u}{\partial t}(t,x) = p\frac{\partial^2 u}{\partial x^2}(t,x) - qu(t,x), \tag{3.1}$$

with nonlocal BCs

$$u(t,0) = u'(t,1) - \int_0^1 \beta_0(s)u(t,s)ds - \int_0^1 \beta_1(s)u'(t,s)ds = 0 \tag{3.2}$$

where $\beta_0, \beta_1 \in L^2(0,1)$ are given functions. Here we are not concerned with existence/uniqueness questions for problem (3.1), (3.2) but rather we assume that $\beta_0, \beta_1 \in L^2(0,1)$ satisfy appropriate conditions so that there exists a set $\tilde{D} \subseteq H^1(0,1)$ with the following property:

"For every $u_0 \in \tilde{D}$ with $u_0(0)=0$, system (3.1), (3.2) has a unique solution $u[t]$ for $t\geq 0$ satisfying $u[0]=u_0$ and $d_1 \in GB$, where $d_1(t) := \int_0^1 \beta_0(s)u(t,s)ds + \int_0^1 \beta_1(s)u'(t,s)ds$, for $t\geq 0$."

Our aim is to provide conditions for exponential stability of system (3.1), (3.2). More specifically, we show that if

$$\|\beta_0\|h(q/p) + \|\beta_1\|\tilde{h}(q/p) < 1 \tag{3.3}$$

where

$$h(s) := \begin{cases} \dfrac{\sqrt{\exp(2\sqrt{s}) - \exp(-2\sqrt{s}) - 4\sqrt{s}}}{(\exp(\sqrt{s}) + \exp(-\sqrt{s}))\sqrt{s}} &, \text{if } s>0 \\ \dfrac{\sqrt{3}}{3} &, \text{if } s=0 \\ \dfrac{\sqrt{2\sqrt{-s} - \sin(2\sqrt{-s})}}{2\cos(\sqrt{-s})(-s)^{3/4}} &, \text{if } -\dfrac{\pi^2}{4} < s < 0 \end{cases}, \quad \tilde{h}(s) := \begin{cases} \dfrac{\sqrt{\exp(2\sqrt{s}) - \exp(-2\sqrt{s}) + 4\sqrt{s}}}{2\sqrt{2}(s)^{1/4}} &, \text{if } s>0 \\ 1 &, \text{if } s=0 \\ \dfrac{\sqrt{2\sqrt{-s} + \sin(2\sqrt{-s})}}{2(-s)^{1/4}\cos(\sqrt{-s})} &, \text{if } -\dfrac{\pi^2}{4} < s < 0 \end{cases} \tag{3.4}$$

then there exist constants $M, \sigma > 0$ such that for every $u_0 \in \tilde{D}$ with $u_0(0)=0$ the unique solution $u[t]$ of system (3.1), (3.2) with initial condition $u[0]=u_0$ satisfies the following decay estimate for all $t\geq 0$:

$$\|u[t]\| + \|u'[t]\| \leq M\exp(-\sigma t)(\|u_0\| + \|u'_0\|) \tag{3.5}$$

In order to prove the decay estimate (3.5), we first study an auxiliary problem. We consider the solution $u[t]$ of (3.1) with

$$u(t,0) = u'(t,1) - d_1(t) = 0 \tag{3.6}$$

The corresponding SL operator $A: D \to L^2(0,1)$ is defined by

$$(Af)(x) = -pf''(x) + qf(x), \text{ for all } f \in D \text{ and } x \in (0,1) \tag{3.7}$$

where $D \subseteq H^2(0,1)$ is the set of all functions $f:[0,1]\to \Re$ for which $f(0)=f'(1)=0$. This operator has eigenfunctions $\phi_n(x) = \sqrt{2}\sin\left(\dfrac{(2n-1)\pi x}{2}\right)$, for $n=1,2,3,...$ and eigenvalues $\lambda_n = \dfrac{(2n-1)^2 \pi^2}{4}p + q$, for $n = 1,2,3,...$. Assumption (H) holds for $A : D \to L^2(0,1)$ and since $q > -p\frac{\pi^2}{4}$, it follows that $\lambda_1 > 0$ and consequently, $A : D \to L^2(0,1)$ is an ES operator. Its derived operator $A' : D' \to L^2(0,1)$ is defined by

$$(A'f)(x) = -p f''(x) + q f(x), \text{ for all } f \in D' \text{ and } x \in (0,1) \tag{3.8}$$

where $D' \subseteq H^2(0,1)$ is the set of all functions $f : [0,1] \to \Re$ for which $f'(0) = f(1) = 0$. This operator has eigenfunctions $\psi_n(x) = \sqrt{2} \cos\left(\frac{(2n-1)\pi x}{2}\right)$, for $n = 1,2,3,...$ and eigenvalues $\mu_n = \frac{(2n-1)^2 \pi^2}{4} p + q$, for $n = 1,2,3,...$. Assumption (H) holds for $A' : D' \to L^2(0,1)$ and since $q > -p\frac{\pi^2}{4}$, it follows that $\mu_1 > 0$ and consequently, $A' : D' \to L^2(0,1)$ is an ES operator. It follows from Theorem 2.9 and Remark 2.10 that for every $d_1 \in GB$, $u_0 \in H^1(0,1)$ with $u_0(0) = 0$, the unique solution $u : \Re_+ \times [0,1] \to \Re$ of the evolution equation (3.1) with (3.6) and initial condition $u_0 \in H^1(0,1)$ satisfies the following estimates for all $\sigma \in \left[0, \frac{\pi^2}{4} p + q\right]$ and $t > 0$:

$$\|u[t]\| \leq \exp(-\sigma t)\|u_0\| + \frac{p\pi^2 + 4q}{p\pi^2 + 4q - 4\sigma} C_1 \sup_{0 < s < t}\left(|d_1(s)|\exp(-\sigma(t-s))\right) \tag{3.9}$$

$$\|u'[t]\| \leq \exp(-\sigma t)\|u'_0\| + \frac{p\pi^2 + 4q}{p\pi^2 + 4q - 4\sigma} \tilde{C}_1 \sup_{0 < s < t}\left(|d_1(s)|\exp(-\sigma(t-s))\right) \tag{3.10}$$

where $C_1 := \|\bar{u}\|$, $\bar{u} \in C^2([0,1])$ is the unique solution of the boundary value problem $p\bar{u}''(x) - q\bar{u}(x) = 0$ for $x \in [0,1]$ with $\bar{u}(0) = 0$ and $\bar{u}'(1) = 1$, $\tilde{C}_1 := \|\bar{v}\|$, $\bar{v} \in C^2([0,1])$ is the solution of the boundary value problem $p\bar{v}''(x) - q\bar{v}(x) = 0$ for $x \in [0,1]$ with $\bar{v}'(0) = 0$ and $\bar{v}(1) = 1$. A direct computation of the solutions of the aforementioned boundary value problems in conjunction with definition (3.4) shows that $C_1 = h(q/p)$ and $\tilde{C}_1 = \tilde{h}(q/p)$. Moreover, if (3.3) holds, then there exists $\sigma > 0$ sufficiently small so that

$$\frac{p\pi^2 + 4q}{p\pi^2 + 4q - 4\sigma} C_1 \|\beta_0\| + \frac{p\pi^2 + 4q}{p\pi^2 + 4q - 4\sigma} \tilde{C}_1 \|\beta_1\| < 1 \tag{3.11}$$

Next, we notice that the solution of (3.1), (3.2) coincides with the solution of (3.1), (3.6) when $d_1(t) = \int_0^1 \beta_0(s) u(t,s) ds + \int_0^1 \beta_1(s) u'(t,s) ds$. Consequently, estimates (3.9), (3.10) hold for all $t > 0$ and for the specific $\sigma > 0$ for which (3.11) holds with $d_1(t) = \int_0^1 \beta_0(s) u(t,s) ds + \int_0^1 \beta_1(s) u'(t,s) ds$. It follows from the Cauchy-Schwarz inequality that

$$|d_1(t)| \leq \|\beta_0\|\|u[t]\| + \|\beta_1\|\|u'[t]\|, \text{ for all } t > 0 \tag{3.12}$$

Combining (3.9), (3.10) and (3.12), we get for all $t > 0$:

$$\|u[t]\|\exp(\sigma t) \leq \|u_0\| + \frac{p\pi^2 + 4q}{p\pi^2 + 4q - 4\sigma} C_1 \|\beta_0\| \sup_{0 < s < t}\left(\|u[s]\|\exp(\sigma s)\right)$$
$$+ \frac{p\pi^2 + 4q}{p\pi^2 + 4q - 4\sigma} C_1 \|\beta_1\| \sup_{0 < s < t}\left(\|u'[s]\|\exp(\sigma s)\right) \tag{3.13}$$

$$\|u'[t]\|\exp(\sigma t) \leq \|u'_0\| + \frac{p\pi^2 + 4q}{p\pi^2 + 4q - 4\sigma} \tilde{C}_1 \|\beta_0\| \sup_{0 < s < t}\left(\|u[s]\|\exp(\sigma s)\right)$$
$$+ \frac{p\pi^2 + 4q}{p\pi^2 + 4q - 4\sigma} \tilde{C}_1 \|\beta_1\| \sup_{0 < s < t}\left(\|u'[s]\|\exp(\sigma s)\right) \tag{3.14}$$

Since (3.11) holds, we get from (3.13) for all $t>0$:

$$\sup_{0<s<t}\left(\|u[s]\|\exp(\sigma s)\right) \leq \left(1 - \frac{p\pi^2+4q}{p\pi^2+4q-4\sigma}C_1\|\beta_0\|\right)^{-1}\left(\|u_0\| + \frac{p\pi^2+4q}{p\pi^2+4q-4\sigma}C_1\|\beta_1\|\sup_{0<s<t}\left(\|u'[s]\|\exp(\sigma s)\right)\right) \quad (3.15)$$

Combining (3.14) and (3.15) and using (3.11), we get for all $t>0$:

$$\sup_{0<s<t}\left(\|u'[s]\|\exp(\sigma s)\right)$$

$$\leq \left(1 - \frac{p\pi^2+4q}{p\pi^2+4q-4\sigma}C_1\|\beta_0\| - \frac{p\pi^2+4q}{p\pi^2+4q-4\sigma}\tilde{C}_1\|\beta_1\|\right)^{-1}\left(1 - \frac{p\pi^2+4q}{p\pi^2+4q-4\sigma}C_1\|\beta_0\|\right)\|u'_0\| \quad (3.16)$$

$$+ \left(1 - \frac{p\pi^2+4q}{p\pi^2+4q-4\sigma}C_1\|\beta_0\| - \frac{p\pi^2+4q}{p\pi^2+4q-4\sigma}\tilde{C}_1\|\beta_1\|\right)^{-1}\frac{p\pi^2+4q}{p\pi^2+4q-4\sigma}\tilde{C}_1\|\beta_0\|\|u_0\|$$

The decay estimate (3.5) with appropriate constant $M>0$ is a direct consequence of estimates (3.15), (3.16). ◁

**Example 3.2:** Let $p>0$, $q\in\Re$ with $q>-p\pi^2$ and consider the nonlocal PDE

$$\frac{\partial u}{\partial t}(t,x) = p\frac{\partial^2 u}{\partial x^2}(t,x) - qu(t,x) + \int_0^1 \beta(s)u'(t,s)ds, \quad (3.17)$$

where $\beta\in L^2(0,1)$ is a given function with Dirichlet BCs

$$u(t,0) = u(t,1) = 0 \quad (3.18)$$

Here, again, we are not concerned with existence/uniqueness questions for problem (3.17), (3.18) but rather we assume that $\beta\in L^2(0,1)$ satisfies appropriate conditions so that there exists a set $\tilde{D}\subseteq H^1(0,1)$ with the following property:

"For every $u_0\in\tilde{D}$ with $u_0(0)=u_0(1)=0$, system (3.17), (3.18) has a unique solution $u[t]$

for all $t\geq 0$ satisfying $u[0]=u_0$, $d\in GB$, where $d(t):=\int_0^1\beta(s)u'(t,s)ds$ for $t\geq 0$."

Notice that if $d\in GB$ then the function $f(t,x):=d(t)$ for $(t,x)\in\Re_+\times[0,1]$ is of class $GD$. Our aim is to provide sufficient conditions for exponential stability of system (3.17), (3.18). More specifically, we show that if

$$p^{-1}h(q/p)\|\beta\| < 1 \quad (3.19)$$

where $h$ is defined by (2.18), then there exist constants $M,\sigma>0$ such that for every $u_0\in\tilde{D}$ with $u_0(0)=u_0(1)=0$ the unique solution $u[t]$ of system (3.17), (3.18) with initial condition $u[0]=u_0$ satisfies the following decay estimate for all $t\geq 0$:

$$\|u[t]\| + \|u'[t]\| \leq M\exp(-\sigma t)\left(\|u_0\| + \|u'_0\|\right) \quad (3.20)$$

In order to prove the decay estimate (3.20), we first study an auxiliary problem. We consider the solution $u[t]$ of

$$\frac{\partial u}{\partial t}(t,x) = p\frac{\partial^2 u}{\partial x^2}(t,x) - qu(t,x) + f(t,x) \quad (3.21)$$

with (3.18), $f\in GD$ with $f[t]\in C^1([0,1])$ for all $t\geq 0$, $f'\in GD$ and $f(t,1)$, $f(t,0)$ being of class $GB$. It follows from Theorem 2.8 that for every $u_0\in H^1(0,1)$ with $u_0(0)=u_0(1)=0$, the unique solution $u:\Re_+\times[0,1]\to\Re$ of the evolution equation (3.21) with (3.18) and initial condition $u_0\in H^1(0,1)$ satisfies estimates (2.16), (2.17) for all $\sigma\in[0,p\pi^2+q)$ and $t>0$. Moreover, if (3.19) holds, then there exists $\sigma>0$ sufficiently small so that

$$\left(\frac{p\pi^2+q}{p\pi^2+q-\sigma}\right)p^{-1}h(q/p)\|\beta\|<1 \tag{3.22}$$

Next, we notice that the solution of (3.17), (3.18) coincides with the solution of (3.21), (3.18) when $f(t,x):=\int_0^1 \beta(s)u'(t,s)ds$. Consequently, estimates (2.16), (2.17) hold for all $t>0$ and for the specific $\sigma>0$ for which (3.22) holds with $f(t,x):=\int_0^1 \beta(s)u'(t,s)ds$. It follows from the Cauchy-Schwarz inequality that

$$\|f[t]\|\leq \|\beta\|\|u'[t]\|, \text{ for all } t\geq 0. \tag{3.23}$$

Combining (2.17) and (3.23), we get for all $t>0$:

$$\|u'[t]\|\exp(\sigma t)\leq \|u'_0\|+\frac{p\pi^2+q}{p\pi^2+q-\sigma}p^{-1}h(q/p)\|\beta\|\sup_{0<s<t}\left(\|u'[s]\|\exp(\sigma s)\right) \tag{3.24}$$

Since (3.22) holds, we get from (3.24) for all $t>0$:

$$\sup_{0<s<t}\left(\|u'[s]\|\exp(\sigma s)\right)\leq \left(1-\frac{p\pi^2+q}{p\pi^2+q-\sigma}p^{-1}h(q/p)\|\beta\|\right)^{-1}\|u'_0\| \tag{3.25}$$

The decay estimate (3.20) with appropriate constant $M>0$ is a direct consequence of estimates (2.16) and (3.25). ◁

# 4. Proofs of Main Results

The proof of Theorem 2.2 requires some technical results. The first two technical results provide a classical solution for a particular parabolic initial-boundary value problem.

**Theorem 4.1:** *Consider the SL operator $A:D\to L^2_r(0,1)$ defined by (2.1), (2.2), where $a_1,a_2,b_1,b_2$ are real constants with $|a_1|+|a_2|>0$, $|b_1|+|b_2|>0$, under Assumption (H). Let $f\in C^0((0,+\infty)\times[0,1])$ be a given function with $\sup_{t\in(0,T)}\left(\int_0^1 r(x)|f(t,x)|^2 dx\right)<+\infty$ for every $T>0$ for which $(0,+\infty)\times[0,1]\ni(t,x)\to\frac{\partial f}{\partial t}(t,x)$ is continuous and $f[t]\in PC^1([0,1])$ for all $t>0$. Then for every $u_0\in L^2_r(0,1)$ and $T>0$, there exists a unique mapping $u\in C^0([0,T];L^2_r(0,1))$, with $u\in C^1((0,T)\times[0,1])$ satisfying $u[t]\in C^2([0,1])$ for all $t\in(0,T]$, $u(0,x)=u_0(x)$ for all $x\in[0,1]$ and*

$$\frac{\partial u}{\partial t}(t,x)-\frac{1}{r(x)}\frac{\partial}{\partial x}\left(p(x)\frac{\partial u}{\partial x}(t,x)\right)+\frac{q(x)}{r(x)}u(t,x)=f(t,x), \text{ for all } (t,x)\in(0,T)\times(0,1) \tag{4.1}$$

$$b_1 u(t,0)+b_2\frac{\partial u}{\partial x}(t,0)=a_1 u(t,1)+a_2\frac{\partial u}{\partial x}(t,1)=0, \text{ for all } t\in(0,T) \tag{4.2}$$

**Proof:** Without loss of generality we assume that $\lambda_1>0$. If this is not the case, we perform the analysis for the function $v(t,x)=\exp(-kt)u(t,x)$ with $k>\lambda_1$ (the function $v(t,x)$ satisfies a PDE similar to (4.1) with the corresponding SL operator satisfying $\lambda_1>0$). It follows from (2.3) that

$$\sum_{n=1}^{\infty}\lambda_n^{-1}\max_{0\leq x\leq 1}\left(|\phi_n(x)|\right)<+\infty \tag{4.3}$$

Define:

$$c_n = \int_0^1 r(x)\phi_n(x)u_0(x)dx, \text{ for } n=1,2,... \quad (4.4)$$

$$\theta_n(t) := \int_0^1 r(x)\phi_n(x)f(t,x)dx, \text{ for all } t>0, \ n=1,2,... \quad (4.5)$$

Since the mapping $[0,1] \ni x \to f(t,x) \in \Re$ is $PC^1([0,1])$ for each $t>0$, it follows from Theorem 11.2.4 in [3], that the following equation holds:

$$f(t,x) = \sum_{n=1}^{\infty} \theta_n(t)\phi_n(x), \text{ for all } (t,x) \in (0,+\infty) \times (0,1) \quad (4.6)$$

Moreover, notice that the Cauchy-Schwarz inequality, in conjunction with the fact that $\|\phi_n\|_r = 1$ (for $n=1,2,...$) and the fact that the mapping $(0,+\infty) \times [0,1] \ni (t,x) \to \frac{\partial f}{\partial t}(t,x)$ is continuous, implies the following relations:

$$|\theta_n(t)| \leq \left( \int_0^1 r(x)|f(t,x)|^2 dx \right)^{1/2}, \text{ for all } t>0 \quad (4.7)$$

$$\dot{\theta}_n(t) = \int_0^1 r(x)\phi_n(x) \frac{\partial f}{\partial t}(t,x)dx, \text{ for } t>0 \quad (4.8)$$

$$|\dot{\theta}_n(t)| \leq \left( \int_0^1 r(x)\left|\frac{\partial f}{\partial t}(t,x)\right|^2 dx \right)^{1/2}, \text{ for } t>0 \quad (4.9)$$

Notice that since the mapping $(0,+\infty) \times [0,1] \ni (t,x) \to \frac{\partial f}{\partial t}(t,x)$ is continuous, it follows that the mappings $\Re_+ \ni t \to \theta_n(t) \in \Re$ is $C^1$ on $(0,+\infty)$. Since for every $T>0$ it holds that $\sup_{t \in (0,T+1)} \left( \int_0^1 r(x)|f(t,x)|^2 dx \right) < +\infty$, it follows that for every $T>0$ there exists $M>0$ such that $\int_0^1 r(x)|f(t,x)|^2 dx \leq M^2$ for all $t \in (0,T]$. It follows from (4.7) that

$$|\theta_n(t)| \leq M, \text{ for all } t \in (0,T] \text{ and } n=1,2,... \quad (4.10)$$

$$\left| \phi_n(x) \int_0^t \exp(-\lambda_n(t-s))\theta_n(s)ds \right| \leq M\lambda_n^{-1} \max_{0 \leq x \leq 1}(|\phi_n(x)|), \text{ for all } (t,x) \in (0,T] \times [0,1] \text{ and } n=1,2,... \quad (4.11)$$

Moreover, the Cauchy-Schwarz inequality and (4.4) imply that $|c_n| \leq \|u_0\|_r$ for $n=1,2,...$ and since $\lambda_n \exp(-\lambda_n t) = t^{-1}\lambda_n t \exp(-\lambda_n t) \leq t_0^{-1} \exp(-1)$ for all $t \in [t_0, T]$ with $t_0 \in (0,T)$, it follows that

$$|\phi_n(z)\exp(-\lambda_n t)c_n| \leq t_0^{-1} \exp(-1)\lambda_n^{-1} \max_{0 \leq x \leq 1}(|\phi_n(x)|)\|u_0\|_r, \text{ for } (t,x) \in [t_0, t_1] \times [0,1], \ n=1,2,... \quad (4.12)$$

Since the Green's function of the SL operator $A: D \to L^2_r(0,1)$ defined by (2.1), (2.2), $g \in C^0([0,1]^2; \Re)$ is a $C^2$ function on each one of the triangles $0 \leq s \leq x \leq 1$ and $0 \leq x \leq s \leq 1$, having a step discontinuity in $\frac{\partial g}{\partial x}(x,s)$ on the line segment $0 \leq s = x \leq 1$ and since $\phi_n(x) = \lambda_n \int_0^1 g(x,s)r(s)\phi_n(s)ds$ for all $x \in [0,1]$, it follows that there exists a constant $K>0$ such that

$$\max_{0 \leq x \leq 1}(|\phi'_n(x)|) \leq K\lambda_n \max_{0 \leq x \leq 1}(|\phi_n(x)|), \text{ for all } n=1,2,... \quad (4.13)$$

The equation $p(x)\phi''_n(x) = q(x)\phi_n(x) - \lambda_n r(x)\phi_n(x) - p'(x)\phi'_n(x)$, which holds for all $x \in [0,1]$, in conjunction with the fact that $0 < \lambda_1 < \lambda_2 < ... < \lambda_n < ...$ and (4.13), implies that that there exists a constant $G>0$ such that

$$\max_{0 \leq x \leq 1}(|\phi''_n(x)|) \leq G\lambda_n \max_{0 \leq x \leq 1}(|\phi_n(x)|), \text{ for all } n=1,2,... \quad (4.14)$$

Moreover, equation (4.6) implies that the following equality holds for all $(t,x) \in (0,+\infty) \times [0,1]$:

$$\sum_{n=1}^{\infty} \phi_n(x) \lambda_n^{-1} \theta_n(t) = \int_0^1 r(s)g(x,s)f(t,s)ds \quad (4.15)$$

Inequalities (4.10), (4.11), (4.12) and (4.3) imply that the series

$$\sum_{n=1}^{\infty} \phi_n(x) \left( \exp(-\lambda_n t) c_n + \int_0^t \exp(-\lambda_n(t-s))\theta_n(s)ds - \lambda_n^{-1}\theta_n(t) \right)$$

is uniformly and absolutely convergent on $[t_0, T] \times [0,1]$ for all $t_0 \in (0,T)$. Therefore, we define $u \in C^0((0,T] \times [0,1])$ by means of the formula:

$$u(t,x) = \int_0^1 g(x,s)r(s)f(t,s)ds + \sum_{n=1}^{\infty} \phi_n(x) \left( \exp(-\lambda_n t) c_n + \int_0^t \exp(-\lambda_n(t-s))\theta_n(s)ds - \lambda_n^{-1}\theta_n(t) \right),$$

for $(t,x) \in (0,T] \times [0,1]$ (4.16)

and we also define

$$u(0,x) := u_0(x), \text{ for all } x \in [0,1] \quad (4.17)$$

The fact that $u_0 \in L_r^2(0,1)$ (which implies that $\sum_{n=1}^{\infty} c_n^2 < +\infty$) in conjunction with (4.10), (4.15), (4.16), (4.17), the fact that $\sum_{n=1}^{\infty} \lambda_n^{-2}$ (a consequence of (4.3) and the fact that

$$1 = \int_0^1 r(x)\phi_n^2(x)dz \leq \left( \max_{0 \leq x \leq 1}(|\phi_n(x)|) \right)^2 \int_0^1 r(x)dx \text{ for } n=1,2,...) \text{ shows that for all } t \in [0,T] \text{ and for every integer}$$

$N \geq 1$, it holds that

$$\|u[t] - u_0\|_r^2 \leq 2\sum_{n=1}^{\infty} (\exp(-\lambda_n t)-1)^2 c_n^2 + 2M^2 \sum_{n=1}^{\infty} \lambda_n^{-2} (\exp(-\lambda_n t)-1)^2$$

$$\leq 2\sum_{n=N+1}^{\infty} c_n^2 + 2M^2 \sum_{n=N+1}^{\infty} \lambda_n^{-2} + 2(\exp(-\lambda_N t)-1)^2 \left( \sum_{n=1}^{N} c_n^2 + M^2 \sum_{n=1}^{\infty} \lambda_n^{-2} \right)$$

The above inequality shows that the mapping $\Re_+ \ni t \to u[t] \in L_r^2([0,1])$ is continuous. We show next that $u[t] \in C^2([0,1])$ for all $t \in (0,T]$ and satisfies the following equations for all $(t,x) \in (0,T] \times [0,1]$:

$$\frac{\partial u}{\partial x}(t,x) = \frac{\partial}{\partial x} \int_0^1 g(x,s)r(s)f(t,s)ds + \sum_{n=1}^{\infty} \phi_n'(x) \left( \exp(-\lambda_n t) c_n + \int_0^t \exp(-\lambda_n(t-s))\theta_n(s)ds - \lambda_n^{-1}\theta_n(t) \right) \quad (4.18)$$

$$\frac{\partial^2 u}{\partial x^2}(t,x) = \frac{\partial^2}{\partial x^2} \int_0^1 g(x,s)r(s)f(t,s)ds + \sum_{n=1}^{\infty} \phi_n''(x) \left( \exp(-\lambda_n t) c_n + \int_0^t \exp(-\lambda_n(t-s))\theta_n(s)ds - \lambda_n^{-1}\theta_n(t) \right) \quad (4.19)$$

This will be achieved by showing that for every $t_0 \in (0, T/2)$ the series appearing in the right hand sides of (4.18) and (4.19) are uniformly and absolutely convergent on $[2t_0, T] \times [0,1]$. Notice that since the mapping $(0,+\infty) \times [0,1] \ni (t,x) \to \frac{\partial f}{\partial t}(t,x)$ is continuous, it follows that for every $t_0 \in (0,T/2)$ there exists $\Gamma > 0$ such that $\int_0^1 r(x) \left| \frac{\partial f}{\partial t}(t,x) \right|^2 dx \leq \Gamma^2$ for all $t \in [t_0, T]$. Therefore, inequality (4.9) implies that

$$\left| \dot{\theta}_n(t) \right| \leq \Gamma, \text{ for all } t \in [t_0, T] \quad (4.20)$$

Using the inequalities $\xi^2 \exp(-\xi) \leq 4\exp(-2)$, $\xi \exp(-\xi) \leq \exp(-1)$ which hold for all $\xi \geq 0$, (4.10), (4.20) and the fact that $|c_n| \leq \|u_0\|_r$ for $n=1,2,...$, we get for all $t \in [2t_0, T]$ and $n=1,2,...$:

$$\left| \exp(-\lambda_n t) c_n + \int_0^t \exp(-\lambda_n(t-s))\theta_n(s)ds - \lambda_n^{-1}\theta_n(t) \right| =$$

$$\left| \exp(-\lambda_n t) c_n + \int_0^{t_0} \exp(-\lambda_n(t-s))\theta_n(s)ds - \lambda_n^{-1}\theta_n(t_0)\exp(-\lambda_n(t-t_0)) - \lambda_n^{-1}\int_{t_0}^t \exp(-\lambda_n(t-s))\dot\theta_n(s)ds \right|$$

$$\leq \exp(-\lambda_n t)|c_n| + \int_0^{t_0} \exp(-\lambda_n(t-s))|\theta_n(s)|ds + \lambda_n^{-1}\exp(-\lambda_n(t-t_0))|\theta_n(t_0)| + \lambda_n^{-1}\int_{t_0}^t \exp(-\lambda_n(t-s))|\dot\theta_n(s)|ds$$

$$\leq 4\lambda_n^{-2} t^{-2}\exp(-2)\|u_0\|_r + Mt_0\exp(-\lambda_n(t-t_0)) + \lambda_n^{-2}(t-t_0)^{-1}\exp(-1)M + \lambda_n^{-2}\Gamma$$

$$\leq \lambda_n^{-2}\left( t_0^{-2}\exp(-2)\|u_0\|_r + 4Mt_0^{-1}\exp(-2) + \exp(-1)Mt_0^{-1} + \Gamma \right)$$

The above inequalities in conjunction with (4.13), (4.14) and (4.3) indicate that for every $t_0 \in (0, T/2)$ the series appearing in the right hand sides of (4.18) and (4.19) are uniformly and absolutely convergent on $[2t_0, T]\times[0,1]$. Similarly, we show that $\dfrac{\partial u}{\partial t}(t,x)$ is a continuous function on $(0,T]\times[0,1]$ and satisfies

$$\frac{\partial u}{\partial t}(t,x) = \int_0^1 g(x,s)r(s)\frac{\partial f}{\partial t}(t,s)ds - \sum_{n=1}^\infty \phi_n(x)\left(\lambda_n\exp(-\lambda_n t)c_n - \theta_n(t) + \lambda_n\int_0^t \exp(-\lambda_n(t-s))\theta_n(s)ds + \lambda_n^{-1}\dot\theta_n(t)\right),$$

$$\text{for all } (t,x) \in (0,T]\times[0,1] \tag{4.21}$$

Again, this is achieved by showing that for every $t_0 \in (0, T/2)$ the series appearing in the right hand side of (4.21) is uniformly and absolutely convergent on $[2t_0, T]\times[0,1]$. Using the inequalities $\xi^2\exp(-\xi) \leq 4\exp(-2)$, $\xi\exp(-\xi)\leq \exp(-1)$ which hold for all $\xi \geq 0$, (4.10), (4.20) and the fact that $|c_n| \leq \|u_0\|_r$ for $n=1,2,...$, we get for all $t \in [2t_0, T]$ and $n = 1,2,...$:

$$\left| \lambda_n\exp(-\lambda_n t)c_n + \lambda_n\int_0^t \exp(-\lambda_n(t-s))\theta_n(s)ds - \theta_n(t) + \lambda_n^{-1}\dot\theta_n(t) \right| =$$

$$\left| \lambda_n\exp(-\lambda_n t)c_n + \lambda_n\int_0^{t_0} \exp(-\lambda_n(t-s))\theta_n(s)ds - \theta_n(t_0)\exp(-\lambda_n(t-t_0)) + \lambda_n^{-1}\dot\theta_n(t) - \int_{t_0}^t \exp(-\lambda_n(t-s))\dot\theta_n(s)ds \right|$$

$$\leq \lambda_n\exp(-\lambda_n t)|c_n| + \lambda_n\int_0^{t_0}\exp(-\lambda_n(t-s))|\theta_n(s)|ds + \exp(-\lambda_n(t-t_0))|\theta_n(t_0)| + \lambda_n^{-1}|\dot\theta_n(t)| + \int_{t_0}^t \exp(-\lambda_n(t-s))|\dot\theta_n(s)|ds$$

$$\leq 4\lambda_n^{-1} t^{-2}\exp(-2)\|u_0\|_r + Mt_0\lambda_n\exp(-\lambda_n(t-t_0)) + \lambda_n^{-1}(t-t_0)^{-1}\exp(-1)M + 2\lambda_n^{-1}\Gamma$$

$$\leq \lambda_n^{-1}\left( t_0^{-2}\exp(-2)\|u_0\|_r + 4Mt_0^{-1}\exp(-2) + \exp(-1)Mt_0^{-1} + 2\Gamma \right)$$

The above inequalities in conjunction with (4.3) indicate that for every $t_0 \in (0, T/2)$ the series appearing in the right hand side of (4.21) is uniformly and absolutely convergent on $[2t_0, T]\times[0,1]$.

Moreover, equalities (4.6), (4.8) imply that the following equality holds for all $(t,x) \in (0,+\infty)\times[0,1]$:

$$\sum_{n=1}^\infty \phi_n(x)\lambda_n^{-1}\dot\theta_n(t) = \int_0^1 r(s)g(x,s)\frac{\partial f}{\partial t}(t,s)ds \tag{4.22}$$

Indeed, for every $t_0 \in (0, T)$ the series appearing in the left hand side of (4.22) is uniformly and absolutely convergent on $[t_0, T]\times[0,1]$ (a consequence of (4.3) and (4.20)).

Equations (4.1), (4.2) are consequences of (4.6), (4.16), (4.21), (4.22), the fact that $A\phi_n = \lambda_n\phi_n$ and the fact that $(y[t])(x) = \int_0^1 r(s)g(x,s)f(t,s)ds$ is the solution of the boundary value problem $(Ay[t])(x) = f(t,x)$ with $b_1(y[t])(0) + b_2\dfrac{d(y[t])}{dx}(0) = a_1(y[t])(1) + a_2\dfrac{d(y[t])}{dx}(1) = 0$ for each fixed $t \geq 0$. Finally, uniqueness follows from Corollary 2.2 on page 106 of the book [42], since $u$ is a strong solution (see Definition 2.8 on page 109 of the above book). The proof is complete. ◁

**Theorem 4.2:** *Consider the SL operator $A: D \to L_r^2(0,1)$ defined by (2.1), (2.2), where $a_1, a_2, b_1, b_2$ are real constants with $|a_1|+|a_2|>0$, $|b_1|+|b_2|>0$, under Assumption (H). Let $f \in C^0((0,+\infty) \times [0,1])$ be a given function with $f[t] \in D \cap C^2([0,1])$ for all $t > 0$ and $\sup_{t \in (0,T)} \left( \int_0^1 r(x) \left( |(Af[t])(x)|^2 + |f(t,x)|^2 \right) dx \right) < +\infty$ for every $T > 0$. Then for every $u_0 \in L_r^2(0,1)$ and $T > 0$, there exists a unique mapping $u \in C^0([0,T]; L_r^2(0,1))$, with $u \in C^1((0,T) \times [0,1])$ satisfying $u[t] \in C^2([0,1])$ for all $t \in (0,T]$, $u(0,x) = u_0(x)$ for all $x \in [0,1]$ and equations (4.1), (4.2).*

**Proof:** Again, without loss of generality, we may replace $A$ with $A + kI$, $k > -\lambda_1$ and assume that $\lambda_1 > 0$ and that $\sup_{t \in (0,T)} \left( \int_0^1 r(x) |(Af[t])(x)|^2 dx \right) < +\infty$ for every $T > 0$. Moreover, we may assume that (4.3) holds. Define $c_n$ for $n = 1, 2, ...$ and $\theta_n(t)$ for $t > 0$, $n = 1, 2, ...$ by (4.4), (4.5), respectively. Notice that since the mapping $[0,1] \ni x \to f(t,x) \in \Re$ is $C^2([0,1])$ for each $t > 0$, Theorem 11.2.4 in [3] implies that the equation (4.6) holds.

Moreover, notice that the Cauchy-Schwarz inequality, in conjunction with the fact that $\|\phi_n\|_r = 1$ (for $n = 1, 2, ...$), definition (4.5) and the fact that $f[t] \in D \cap C^2([0,1])$ for all $t > 0$, implies that

$$|\theta_n(t)| \leq \lambda_n^{-1} \left( \int_0^1 r(x) |(Af[t])(x)|^2 dx \right)^{1/2}, \text{ for all } t > 0 \quad (4.23)$$

Notice that since $f \in C^0((0,+\infty) \times [0,1])$, it follows that the mappings $\Re_+ \ni t \to \theta_n(t) \in \Re$ is $C^0$ on $(0,+\infty)$. Since for every $T > 0$ it holds that $\sup_{t \in (0,T+1)} \left( \int_0^1 r(x) |(Af[t])(x)|^2 dx \right) < +\infty$, it follows that for every $T > 0$ there exists $M > 0$ such that $\int_0^1 r(x) |(Af[t])(x)|^2 dx \leq M^2$ for all $t \in (0,T]$. It follows from (4.23) that

$$|\theta_n(t)| \leq \lambda_n^{-1} M, \text{ for all } t \in (0,T] \text{ and } n = 1, 2, ... \quad (4.24)$$

$$\left| \phi_n(x) \int_0^t \exp(-\lambda_n(t-s)) \theta_n(s) ds \right| \leq M \lambda_n^{-2} \max_{0 \leq x \leq 1}(|\phi_n(x)|), \text{ for all } (t,x) \in (0,T] \times [0,1] \text{ and } n = 1, 2, ... \quad (4.25)$$

Moreover, the Cauchy-Schwarz inequality and (4.4) imply that $|c_n| \leq \|u_0\|_r$ for $n = 1, 2, ...$ and since $\lambda_n \exp(-\lambda_n t) = t^{-1} \lambda_n t \exp(-\lambda_n t) \leq t_0^{-1} \exp(-1)$ for all $t \in [t_0, T]$ with $t_0 \in (0,T)$, it follows that (4.12) holds.

Inequalities (4.25), (4.12) and (4.3) imply that the series

$$\sum_{n=1}^\infty \phi_n(x) \left( \exp(-\lambda_n t) c_n + \int_0^t \exp(-\lambda_n(t-s)) \theta_n(s) ds \right)$$

is uniformly and absolutely convergent on $[t_0, T] \times [0,1]$ for all $t_0 \in (0,T)$. Therefore, we define $u \in C^0((0,T] \times [0,1])$ by means of the formula:

$$u(t,x) = \sum_{n=1}^\infty \phi_n(x) \left( \exp(-\lambda_n t) c_n + \int_0^t \exp(-\lambda_n(t-s)) \theta_n(s) ds \right),$$

$$\text{for } (t,x) \in (0,T] \times [0,1] \quad (4.26)$$

and we also define

$$u(0,x) := u_0(x), \text{ for all } x \in [0,1] \quad (4.27)$$

The fact that $u_0 \in L^2_r(0,1)$ (which implies that $\sum_{n=1}^{\infty} c_n^2 < +\infty$) in conjunction with (4.26), the fact that $\sum_{n=1}^{\infty} \lambda_n^{-2}$ (a consequence of (4.3) and the fact that $1 = \int_0^1 r(x)\phi_n^2(x)dx \leq \left(\max_{0\leq x\leq 1}(|\phi_n(x)|)\right)^2 \int_0^1 r(x)dx$ for $n=1,2,...$) shows that for all $t \in [0,T]$ and for every integer $N \geq 1$, it holds that

$$\|u[t] - u_0\|_r^2 \leq 2\sum_{n=1}^{\infty}(\exp(-\lambda_n t)-1)^2 c_n^2 + 2M^2 \sum_{n=1}^{\infty} \lambda_n^{-2}(\exp(-\lambda_n t)-1)^2$$

$$\leq 2\sum_{n=N+1}^{\infty} c_n^2 + 2M^2 \sum_{n=N+1}^{\infty} \lambda_n^{-2} + 2(\exp(-\lambda_N t)-1)^2 \left(\sum_{n=1}^{N} c_n^2 + M^2 \sum_{n=1}^{\infty} \lambda_n^{-2}\right)$$

The above inequality shows that the mapping $\Re_+ \ni t \to u[t] \in L^2_r([0,1])$ is continuous.

We show next that $u[t] \in C^2([0,1])$ for all $t \in (0,T]$ and satisfies the following equations for all $(t,x) \in (0,T] \times [0,1]$:

$$\frac{\partial u}{\partial x}(t,x) = \sum_{n=1}^{\infty} \phi_n'(x)\left(\exp(-\lambda_n t)c_n + \int_0^t \exp(-\lambda_n(t-s))\theta_n(s)ds\right) \quad (4.28)$$

$$\frac{\partial^2 u}{\partial x^2}(t,x) = \sum_{n=1}^{\infty} \phi_n''(x)\left(\exp(-\lambda_n t)c_n + \int_0^t \exp(-\lambda_n(t-s))\theta_n(s)ds\right) \quad (4.29)$$

This is achieved by showing that for every $t_0 \in (0, T/2)$ the series appearing in the right hand sides of (4.28) and (4.29) are uniformly and absolutely convergent on $[2t_0, T] \times [0,1]$.

Indeed, there exist constants $K, G > 0$ such that inequalities (4.13), (4.14) hold for all $n=1,2,...$. Using the inequality $\xi^2 \exp(-\xi) \leq 4\exp(-2)$, which holds for all $\xi \geq 0$, (4.24) and the fact that $|c_n| \leq \|u_0\|_r$ for $n=1,2,...$, we get for all $t \in [2t_0, T]$ and $n=1,2,...$:

$$\left|\exp(-\lambda_n t)c_n + \int_0^t \exp(-\lambda_n(t-s))\theta_n(s)ds\right| \leq \exp(-\lambda_n t)|c_n| + \int_0^t \exp(-\lambda_n(t-s))|\theta_n(s)|ds$$

$$\leq 4\lambda_n^{-2}t^{-2}\exp(-2)\|u_0\|_r + \lambda_n^{-2}M \leq \lambda_n^{-2}\left(t_0^{-2}\exp(-2)\|u_0\|_r + M\right)$$

The above inequalities in conjunction with (4.13), (4.14) and (4.3) indicate that for every $t_0 \in (0, T/2)$ the series appearing in the right hand sides of (4.28) and (4.29) are uniformly and absolutely convergent on $[2t_0, T] \times [0,1]$. Similarly, we show that $\frac{\partial u}{\partial t}(t,x)$ is a continuous function on $(0,T] \times [0,1]$ and satisfies

$$\frac{\partial u}{\partial t}(t,x) = -\sum_{n=1}^{\infty} \phi_n(x)\left(\lambda_n \exp(-\lambda_n t)c_n - \theta_n(t) + \lambda_n \int_0^t \exp(-\lambda_n(t-s))\theta_n(s)ds\right), \text{ for all } (t,x) \in (0,T] \times [0,1] \quad (4.30)$$

Equations (4.1), (4.2) are consequences of (4.6), (4.26), (4.30) and the fact that $A\phi_n = \lambda_n \phi_n$. Finally, uniqueness follows from Corollary 2.2 on page 106 of the book [42], since the constructed function $u$ is a strong solution (see Definition 2.8 on page 109 of the above book).
The proof is complete. ◁

The following technical result allows the construction of a solution for the evolution equation (2.4) with (2.5). Theorem 2.2 follows as a direct consequence of Theorem 4.3, since we can repeat the construction on every interval $(\tau_{i-1}, \tau_i)$ with $i \geq 1$.

**Theorem 4.3:** *Consider the SL operator $A: D \to L^2_r(0,1)$ defined by (2.1), (2.2), where $a_1, a_2, b_1, b_2$ are real constants with $|a_1| + |a_2| > 0$, $|b_1| + |b_2| > 0$, under Assumption (H). Let $T > 0$ be a constant and* let $f_1 \in C^0((0,T)\times[0,1])$, $a \in PC^1([0,1])$, $b \in C^0([0,1])$ be given functions with $\sup_{t\in(0,T)}\left(\int_0^1 r(x)|f_1(t,x)|^2 dx\right) < +\infty$, for which $(0,T)\times[0,1] \ni (t,x) \to \frac{\partial f_1}{\partial t}(t,x)$ is continuous, $f_1[t] \in PC^1([0,1])$ for all $t \in (0,T)$ and $\lim_{t \to T^-}(f_1(t,x)) = a(x)$, $\lim_{t \to T^-}\left(\frac{\partial f_1}{\partial t}(t,x)\right) = b(x)$ for $x \in [0,1]$. Let $f_2 \in C^0((0,T)\times[0,1])$, $c \in D \cap C^2([0,1])$ be given functions with $f_2[t] \in D \cap C^2([0,1])$ for all $t \in (0,T)$, $\sup_{t\in(0,T)}\left(\int_0^1 r(x)\left(|(Af_2[t])(x)|^2 + |f_2(t,x)|^2\right)dx\right) < +\infty$ and $\lim_{t \to T^-}(f_2(t,x)) = c(x)$. Let $d_0, d_1 : [0,T) \to \Re$ be right continuous functions, $C^2$ functions on $(0,T)$ with $\sup_{t\in(0,T)}\left(|\dot{d}_0(t)| + |\dot{d}_1(t)|\right) < +\infty$ and all left limits of $d_0(t), \dot{d}_0(t), \ddot{d}_0(t)$, $d_1(t), \dot{d}_1(t), \ddot{d}_1(t)$ when $t$ tends $T$ being finite. Define $f(t,x) := f_1(t,x) + f_2(t,x)$ for all $(t,x) \in (0,T)\times[0,1]$. Then for every $u_0 \in L_r^2(0,1)$, there exists a unique function $u:[0,T]\times[0,1]\to\Re$ for which the mapping $[0,T] \ni t \to u[t] \in L_r^2(0,1)$ is continuous, with $u \in C^1((0,T)\times[0,1])$ satisfying $u[t] \in C^2([0,1])$ for all $t \in (0,T]$, $\lim_{t \to T^-}(u(t,x)) = u(T,x)$, $\lim_{t \to T^-}\left(\frac{\partial u}{\partial t}(t,x)\right) = -(Au[T])(x) + a(x) + c(x)$, $\lim_{t \to T^-}\left(\frac{\partial u}{\partial x}(t,x)\right) = \frac{\partial u}{\partial x}(T,x)$, $u(0,x) = u_0(x)$ for all $x \in [0,1]$ and

$$\frac{\partial u}{\partial t}(t,x) - \frac{1}{r(x)}\frac{\partial}{\partial x}\left(p(x)\frac{\partial u}{\partial x}(t,x)\right) + \frac{q(x)}{r(x)}u(t,x) = f(t,x), \text{ for all } (t,x) \in (0,T)\times(0,1) \quad (4.31)$$

$$b_1 u(t,0) + b_2 \frac{\partial u}{\partial x}(t,0) - d_0(t) = a_1 u(t,1) + a_2 \frac{\partial u}{\partial x}(t,1) - d_1(t) = 0, \text{ for all } t \in (0,T) \quad (4.32)$$

**Proof:** First, we extend $d_0, d_1$ on $\Re_+$ so that $d_0, d_1 \in C^2((0,+\infty))$ (e.g., $d_0(t) = \lim_{l \to T^-}(d_0(l)) + (t-T)\lim_{l \to T^-}(\dot{d}_0(l)) + \frac{1}{2}(t-T)^2 \lim_{l \to T^-}(\ddot{d}_0(l))$, for $t \in [T,+\infty)$ and similarly for $d_1$). Notice that right continuity of $d_0, d_1 : [0,T) \to \Re$ implies continuity of $d_0, d_1$ on $\Re_+$. Next, we extend $f_1$ to $(0,+\infty)\times[0,1]$ so that $f_1 \in C^0((0,+\infty)\times[0,1])$ with $\sup_{l\in(0,t)}\left(\int_0^1 r(x)|f_1(l,x)|^2 dx\right) < +\infty$ for every $t > 0$, $(0,+\infty)\times[0,1] \ni (t,x) \to \frac{\partial f_1}{\partial t}(t,x)$ being a continuous mapping and $f_1[t] \in PC^1([0,1])$ for all $t > 0$. This can be done by first extending $b$ on $\Re_+$ continuously (e.g., $b(x) = b(1)$ for $x > 1$) and setting

$$f_1(t,x) = a(x) + \int_x^{x+t-T} b(s)ds, \text{ for } (t,x) \in [T,+\infty)\times[0,1].$$

Next, we extend $f_2$ to $(0,+\infty)\times[0,1]$ so that $f_2 \in C^0((0,+\infty)\times[0,1])$ with $f_2[t] \in D \cap C^2([0,1])$ and $\sup_{l\in(0,t)}\left(\int_0^1 r(x)\left(|(Af[l])(x)|^2 + |f(l,x)|^2\right)dx\right) < +\infty$ for all $t > 0$ (e.g., by setting $f_2(t,x) = c(x)$, for $(t,x) \in [T,+\infty)\times[0,1]$).

Finally, we extend $f$ to $(0,+\infty)\times[0,1]$ by setting $f(t,x) := f_1(t,x) + f_2(t,x)$ for all $(t,x) \in (0,+\infty)\times[0,1]$.

Define the function

$$\tilde{f}(t,x) := -\dot{d}_0(t)p_0(x) - \dot{d}_1(t)p_1(x) - d_0(t)(Ap_0)(x) - d_1(t)(Ap_1)(x) + f(t,x),$$

for all $(t,x) \in (0,+\infty)\times[0,1]$, $\quad (4.33)$

where $p_0(x) := (c_0 + c_1 x)(1-x)^2$, $p_1(z) := (C_0 + C_1(1-x))x^2$, $c_0 = \frac{b_1 - 2b_2}{(b_1 - 2b_2)^2 + b_2^2}$, $c_1 = \frac{b_2}{(b_1 - 2b_2)^2 + b_2^2}$, $C_0 = \frac{a_1 + 2a_2}{(a_1 + 2a_2)^2 + a_2^2}$, $C_1 = \frac{-a_2}{(a_1 + 2a_2)^2 + a_2^2}$ (notice that the fact $|a_1| + |a_2| > 0$, $|b_1| + |b_2| > 0$ guarantees that $(b_1 - 2b_2)^2 + b_2^2 > 0$, $(a_1 + 2a_2)^2 + a_2^2 > 0$). The facts that $\sup_{t \in (0,T)} \left( |d_0(t)| + |\dot{d}_0(t)| + |d_1(t)| + |\dot{d}_1(t)| \right) < +\infty$, $\sup_{l \in (0,t)} \left( \int_0^1 r(x) |f(l,x)|^2 dx \right) < +\infty$ for every $t > 0$ and definition (4.33) guarantee that $\tilde{f} \in C^0((0,+\infty) \times [0,1])$ is a function that can be written as the sum of two functions:

- a function $g_1 \in C^0((0,+\infty) \times [0,1])$ with $\sup_{l \in (0,t)} \left( \int_0^1 r(x) |g_1(l,x)|^2 dx \right) < +\infty$ for every $t > 0$ for which $(0,+\infty) \times [0,1] \ni (t,x) \to \frac{\partial g_1}{\partial t}(t,x)$ is continuous and $g_1[t] \in PC^1([0,1])$ for all $t > 0$,

- a function $g_2 \in C^0((0,+\infty) \times [0,1])$ with $g_2[t] \in D \cap C^2([0,1])$ and $\sup_{l \in (0,t)} \left( \int_0^1 r(x) \left( |(Ag_2[l])(x)|^2 + |g_2(l,x)|^2 \right) dx \right) < +\infty$ for all $t > 0$.

Theorem 4.1 and Theorem 4.2 guarantee that for every $u_0 \in L_r^2(0,1)$ and $\tilde{T} > T$, there exists a unique function $y : [0,\tilde{T}] \times [0,1] \to \Re$ for which the mapping $[0,T] \ni t \to y[t] \in L_r^2(0,1)$ is continuous, with $y \in C^1((0,\tilde{T}) \times [0,1])$ satisfying $y[t] \in C^2([0,1])$ for all $t \in (0,\tilde{T}]$, $y(0,x) = u(0,x) - p_0(x)d_0(0) - p_1(x)d_1(0)$ for all $x \in [0,1]$ and

$$\frac{\partial y}{\partial t}(t,x) + (Ay[t])(x) = \tilde{f}(t,x), \text{ for all } (t,x) \in (0,\tilde{T}) \times (0,1) \tag{4.34}$$

$$b_1 y(t,0) + b_2 \frac{\partial y}{\partial x}(t,0) = a_1 y(t,1) + a_2 \frac{\partial y}{\partial x}(t,1) = 0, \text{ for all } t \in (0,\tilde{T}) \tag{4.35}$$

Finally, we define
$$u(t,x) = y(t,x) + p_0(x)d_0(t) + d_1(t)p_1(x), \text{ for } (t,x) \in [0,T] \times [0,1] \tag{4.36}$$

Notice that continuity of $d_0, d_1$ on $\Re_+$ implies that the mapping $[0,T] \ni t \to u[t] \in L_r^2(0,1)$ is continuous. Moreover, (4.34) and definitions (4.33), (4.36) guarantee that $\frac{\partial u}{\partial t}(t,x) = -(Au[t])(x) + f(t,x)$ for all $(t,x) \in (0,T] \times (0,1)$. Continuity with respect to $x$ of $-(Au[t])(x) + f(t,x)$ and $\frac{\partial u}{\partial t}(t,x)$ for $t \in (0,T]$ imply that $\frac{\partial u}{\partial t}(t,x) = -(Au[t])(x) + f(t,x)$ for all $(t,x) \in (0,T] \times [0,1]$. The fact that $f(T,x) = a(x) + c(x)$ implies that $\lim_{t \to T^-} \left( \frac{\partial u}{\partial t}(t,x) \right) = -(Au[T])(x) + a(x) + c(x)$ for all $x \in [0,1]$. All rest conclusions are consequences of (4.35) and definitions (4.33), (4.36). The proof is complete. ◁

We next continue with the proofs of Theorem 2.4 and Proposition 2.6.

**Proof of Theorem 2.4:** First, we prove the following Claim.

**Claim:** For every solution $u : \Re_+ \times [0,1] \to \Re$ of the evolution equation (2.4) with (2.5), initial condition $u_0 \in L_r^2(0,1)$ and corresponding to inputs $f \in GD$, $d_0, d_1 \in GB$, the following equations hold for all $t \geq 0$ and $n = 1,2,...$:

$$c_n(t) = \exp(-\lambda_n t) c_n(0) + \frac{p(1)}{a_1^2 + a_2^2} \left( a_2 \phi_n(1) - a_1 \frac{d\phi_n}{dz}(1) \right) \int_0^t \exp(-\lambda_n(t-s)) d_1(s) ds$$
$$+ \frac{p(0)}{b_1^2 + b_2^2} \left( b_1 \frac{d\phi_n}{dz}(0) - b_2 \phi_n(0) \right) \int_0^t \exp(-\lambda_n(t-s)) d_0(s) ds + \int_0^t \exp(-\lambda_n(t-s)) \left( \int_0^1 r(x) f(s,x) \phi_n(x) dx \right) ds \tag{4.37}$$

where

$$c_n(t) := \int_0^1 r(x)u(t,x)\phi_n(x)dx, \text{ for } n=1,2,... \qquad (4.38)$$

**Proof of Claim:** Let $\{\tau_i \geq 0, i=0,1,2,...\}$ be the increasing sequence of times with $\tau_0 = 0$, $\lim_{i\to+\infty}(\tau_i) = +\infty$ involved in Definition 2.1. Theorem 2.2 guarantees that the mappings $\Re_+ \ni t \to c_n(t)$ are continuous for $n=1,2,...$. Moreover, Theorem 2.2 guarantees that the mappings $I \ni t \to c_n(t)$, where $I = \Re_+ \setminus \{\tau_i \geq 0, i=0,1,2,...\}$, are $C^1$ on $I$. By virtue of (2.6), it follows from repeated integration by parts, that the following equalities hold for all $t \in I$ and $n=1,2,...$:

$$\dot{c}_n(t) = \int_0^1 r(x)\frac{\partial u}{\partial t}(t,x)\phi_n(x)dx$$

$$= p(1)\left(\frac{\partial u}{\partial x}(t,1)\phi_n(1) - x(t,1)\frac{d\phi_n}{dx}(1)\right) + p(0)\left(\frac{d\phi_n}{dz}(0)u(t,0) - \frac{\partial u}{\partial x}(t,0)\phi_n(0)\right)$$

$$- \int_0^1 r(x)u(t,x)(A\phi_n)(x)dx + \int_0^1 r(x)f(t,x)\phi_n(x)dx$$

Thus we get for all $t \in I$ and $n=1,2,...$:

$$\dot{c}_n(t) = p(1)\left(\frac{\partial u}{\partial x}(t,1)\phi_n(1) - u(t,1)\frac{d\phi_n}{dx}(1)\right) + p(0)\left(\frac{d\phi_n}{dx}(0)u(t,0) - \frac{\partial u}{\partial x}(t,0)\phi_n(0)\right)$$
$$- \int_0^1 r(x)u(t,x)(A\phi_n)(x)dx + \int_0^1 r(x)f(t,x)\phi_n(x)dx \qquad (4.39)$$

It follows from (4.39), the fact that $(A\phi_n)(x) = \lambda_n \phi_n(x)$ and definition (4.38) that the following equation holds for all $t \in I$ and $n=1,2,...$:

$$\dot{c}_n(t) + \lambda_n c_n(t) = p(1)\left(\frac{\partial u}{\partial x}(t,1)\phi_n(1) - u(t,1)\frac{d\phi_n}{dx}(1)\right)$$
$$+ p(0)\left(\frac{d\phi_n}{dx}(0)u(t,0) - \frac{\partial u}{\partial x}(t,0)\phi_n(0)\right) + \int_0^1 r(x)f(t,x)\phi_n(x)dx \qquad (4.40)$$

Proceeding exactly as in the proof of Theorem 2.3 in [23], it can be shown that for all $t \in I$ and $n=1,2,...$ the following equalities hold:

$$u(t,0)\frac{d\phi_n}{dx}(0) - \frac{\partial u}{\partial x}(t,0)\phi_n(0) = \frac{d_0(t)}{b_1^2 + b_2^2}\left(b_1\frac{d\phi_n}{dx}(0) - b_2\phi_n(0)\right) \qquad (4.41)$$

$$\phi_n(1)\frac{\partial u}{\partial x}(t,1) - u(t,1)\frac{d\phi_n}{dx}(1) = \frac{d_1(t)}{a_1^2 + a_2^2}\left(a_2\phi_n(1) - a_1\frac{d\phi_n}{dx}(1)\right) \qquad (4.42)$$

Using (4.40), (4.41) and (4.41), we obtain for all $t \in I$ and $n=1,2,...$:

$$\dot{c}_n(t) + \lambda_n c_n(t) = \frac{p(1)d_1(t)}{a_1^2 + a_2^2}\left(a_2\phi_n(1) - a_1\frac{d\phi_n}{dx}(1)\right)$$
$$+ \frac{p(0)d_0(t)}{b_1^2 + b_2^2}\left(b_1\frac{d\phi_n}{dx}(0) - b_2\phi_n(0)\right) + \int_0^1 r(x)f(t,x)\phi_n(x)dx \qquad (4.43)$$

Integrating the differential equations (4.43), we obtain (4.37) for all $t \geq 0$ and $n=1,2,...$. The proof of the claim is complete. ◁

Next, we recognize that the solution $u: \Re_+ \times [0,1] \to \Re$ of the evolution equation (2.4) with (2.5), initial condition $u_0 \in L_r^2(0,1)$ and corresponding to inputs $f \in GD$, $d_0, d_1 \in GB$ satisfies the following equation:

$$u[t] = u_1[t] + u_2[t] + u_3[t] + u_4[t], \text{ for all } t \geq 0 \qquad (4.44)$$

where

- $u_1 : \Re_+ \times [0,1] \to \Re$ is the solution of the evolution equation (2.4) with (2.5), initial condition $u_0 \in L_r^2(0,1)$ and corresponding to inputs $f \equiv 0$, $d_0 = d_1 \equiv 0$,
- $u_2 : \Re_+ \times [0,1] \to \Re$ is the solution of the evolution equation (2.4) with (2.5), initial condition $u_0 \equiv 0$ and corresponding to inputs $f \equiv 0$, $d_0 \in GB$, $d_1 \equiv 0$,
- $u_3 : \Re_+ \times [0,1] \to \Re$ is the solution of the evolution equation (2.4) with (2.5), initial condition $u_0 \equiv 0$ and corresponding to inputs $f \equiv 0$, $d_1 \in GB$, $d_0 \equiv 0$,
- $u_4 : \Re_+ \times [0,1] \to \Re$ is the solution of the evolution equation (2.4) with (2.5), initial condition $u_0 \equiv 0$ and corresponding to inputs $f \in GD$, $d_0 = d_1 \equiv 0$.

Next, we estimate each component of $u : \Re_+ \times [0,1] \to \Re$ separately.

1) Estimate for $u_1 : \Re_+ \times [0,1] \to \Re$.

We obtain from (4.37) for all $t \geq 0$ and $n = 1,2,\ldots$:

$$c_n(t) = \exp(-\lambda_n t) c_n(0) \tag{4.45}$$

where $c_n(t) := \int_0^1 r(x) u_1(t,x) \phi_n(x) dx$. Since the eigenfunctions $\{\phi_n\}_{n=1}^{\infty}$ of the SL operator $A : D \to L_r^2(0,1)$ defined by (2.1), (2.2) form an orthonormal basis of $L_r^2(0,1)$, it follows that Parseval's identity holds, i.e.,

$$\|u_1[t]\|_r^2 = \sum_{n=1}^{\infty} c_n^2(t), \text{ for all } t \geq 0 \tag{4.46}$$

Using (4.45), (4.46) and the fact that $\lambda_n \geq \lambda_1$ for all $n = 1,2,\ldots$, we get for all $t \geq 0$:

$$\|u_1[t]\|_r \leq \exp(-\lambda_1 t) \|u_0\|_r \tag{4.47}$$

2) Estimate for $u_2 : \Re_+ \times [0,1] \to \Re$.

We obtain from (4.37) for all $t \geq 0$ and $n = 1,2,\ldots$:

$$c_n(t) = \frac{p(0)}{b_1^2 + b_2^2} \left( b_1 \frac{d\phi_n}{dx}(0) - b_2 \phi_n(0) \right) \int_0^t \exp(-\lambda_n (t-s)) d_0(s) ds \tag{4.48}$$

where $c_n(t) := \int_0^1 r(x) u_2(t,x) \phi_n(x) dx$. Using (4.48), we get all $\sigma \in [0, \lambda_1)$, $t > 0$ and $n = 1,2,\ldots$:

$$|c_n(t)|^2 \leq \frac{1}{(\lambda_n - \sigma)^2} \left( \frac{p(0)}{b_1^2 + b_2^2} \right)^2 \left| b_1 \frac{d\phi_n}{dx}(0) - b_2 \phi_n(0) \right|^2 \sup_{0<s<t} \left( |d_0(s)|^2 \exp(-2\sigma(t-s)) \right) \tag{4.49}$$

Since the eigenfunctions $\{\phi_n\}_{n=1}^{\infty}$ of the SL operator $A : D \to L_r^2(0,1)$ defined by (2.1), (2.2) form an orthonormal basis of $L_r^2(0,1)$, it follows that Parseval's identity holds, i.e.,

$$\|u_2[t]\|_r^2 = \sum_{n=1}^{\infty} c_n^2(t), \text{ for all } t \geq 0 \tag{4.50}$$

Therefore, by virtue of (4.49), (4.50), the following estimate holds for all $\sigma \in [0, \lambda_1)$, $t > 0$:

$$\|u_2[t]\|_r \leq K_0 \sup_{0<s<t} \left( |d_0(s)| \exp(-\sigma(t-s)) \right) \tag{4.51}$$

where

$$K_0 := \frac{p(0)}{\sqrt{b_1^2 + b_2^2}} \sqrt{ \sum_{n=1}^{\infty} \frac{1}{(\lambda_n - \sigma)^2} \left| \frac{b_1}{\sqrt{b_1^2 + b_2^2}} \frac{d\phi_n}{dx}(0) - \frac{b_2}{\sqrt{b_1^2 + b_2^2}} \phi_n(0) \right|^2 } \tag{4.52}$$

Equation (2.9) is a direct consequence of Theorem 2.3 in [23]. Moreover, definitions (2.9), (4.52) imply that $K_0 \leq \frac{\lambda_1}{\lambda_1 - \sigma} C_0$. Consequently, the following estimate holds for all $\sigma \in [0, \lambda_1)$, $t > 0$:

$$\|u_2[t]\|_r \leq \frac{\lambda_1}{\lambda_1 - \sigma} C_0 \sup_{0 < s < t} \left( |d_0(s)| \exp(-\sigma(t-s)) \right) \tag{4.53}$$

3) Estimate for $u_3 : \mathfrak{R}_+ \times [0,1] \to \mathfrak{R}$.

We obtain from (4.37) for all $t \geq 0$ and $n = 1, 2, \ldots$:

$$c_n(t) = \frac{p(1)}{a_1^2 + a_2^2} \left( a_2 \phi_n(1) - a_1 \frac{d\phi_n}{dx}(1) \right) \int_0^t \exp(-\lambda_n(t-s)) d_1(s) ds \tag{4.54}$$

where $c_n(t) := \int_0^1 r(x) u_3(t,x) \phi_n(x) dx$. Using (4.54), we get all $\sigma \in [0, \lambda_1)$, $t > 0$ and $n = 1, 2, \ldots$:

$$|c_n(t)|^2 \leq \frac{1}{(\lambda_n - \sigma)^2} \left( \frac{p(1)}{a_1^2 + a_2^2} \right)^2 \left| a_1 \frac{d\phi_n}{dx}(0) - a_2 \phi_n(0) \right|^2 \sup_{0<s<t} \left( |d_1(s)|^2 \exp(-2\sigma(t-s)) \right) \tag{4.55}$$

Since the eigenfunctions $\{\phi_n\}_{n=1}^\infty$ of the SL operator $A : D \to L_r^2(0,1)$ defined by (2.1), (2.2) form an orthonormal basis of $L_r^2(0,1)$, it follows that Parseval's identity holds, i.e.,

$$\|u_3[t]\|_r^2 = \sum_{n=1}^\infty c_n^2(t), \text{ for all } t \geq 0 \tag{4.56}$$

Therefore, by virtue of (4.55), (4.56), the following estimate holds for all $\sigma \in [0, \lambda_1)$, $t > 0$:

$$\|u_3[t]\|_r \leq K_1 \sup_{0<s<t} \left( |d_1(s)| \exp(-\sigma(t-s)) \right) \tag{4.57}$$

where

$$K_1 := \frac{p(1)}{\sqrt{a_1^2 + a_2^2}} \sqrt{\sum_{n=1}^\infty \frac{1}{(\lambda_n - \sigma)^2} \left| \frac{a_2}{\sqrt{a_1^2 + a_2^2}} \phi_n(1) - \frac{a_1}{\sqrt{a_1^2 + a_2^2}} \frac{d\phi_n}{dx}(1) \right|^2} \tag{4.58}$$

Equation (2.10) is a direct consequence of Theorem 2.3 in [23]. Moreover, definitions (2.10), (4.58) imply that $K_1 \leq \frac{\lambda_1}{\lambda_1 - \sigma} C_1$. Consequently, the following estimate holds for all $\sigma \in [0, \lambda_1)$, $t > 0$:

$$\|u_3[t]\|_r \leq \frac{\lambda_1}{\lambda_1 - \sigma} C_1 \sup_{0<s<t} \left( |d_1(s)| \exp(-\sigma(t-s)) \right) \tag{4.59}$$

4) Estimate for $u_4 : \mathfrak{R}_+ \times [0,1] \to \mathfrak{R}$.

We obtain from (4.37) for all $t \geq 0$ and $n = 1, 2, \ldots$:

$$c_n(t) = \int_0^t \exp(-\lambda_n(t-s)) f_n(s) ds \tag{4.60}$$

where $c_n(t) := \int_0^1 r(x) u_4(t,x) \phi_n(x) dx$ and $f_n(t) := \int_0^1 r(x) f(t,x) \phi_n(x) dx$. Since the eigenfunctions $\{\phi_n\}_{n=1}^\infty$ of the SL operator $A : D \to L_r^2(0,1)$ defined by (2.1), (2.2) form an orthonormal basis of $L_r^2(0,1)$, it follows that Parseval's identity holds, i.e.,

$$\|u_4[t]\|_r^2 = \sum_{n=1}^\infty c_n^2(t), \text{ for all } t \geq 0 \tag{4.61}$$

$$\|f[t]\|_r^2 = \sum_{n=1}^{\infty} f_n^2(t), \text{ for all } t > 0 \qquad (4.62)$$

It follows from (4.60) that the following inequality holds for all $\sigma \in [0, \lambda_1)$, $t \geq 0$ and $n = 1, 2, \ldots$:

$$|c_n(t)| \leq \int_0^t \exp(-\lambda_n (t-s))|f_n(s)|ds \leq \int_0^t \exp(-(\lambda_1 - \sigma)(t-s))|f_n(s)| \exp(-\sigma(t-s))ds \qquad (4.63)$$

Using the Cauchy-Schwarz inequality and (4.63) we obtain for all $\sigma \in [0, \lambda_1)$, $t \geq 0$ and $n = 1, 2, \ldots$:

$$|c_n(t)| \leq \frac{1}{\sqrt{\lambda_1 - \sigma}} \left( \int_0^t \exp(-(\lambda_1 - \sigma)(t-s))|f_n(s)|^2 \exp(-2\sigma(t-s))ds \right)^{1/2}$$

which directly implies

$$|c_n(t)|^2 \leq \frac{1}{\lambda_1 - \sigma} \int_0^t \exp(-(\lambda_1 - \sigma)(t-s))|f_n(s)|^2 \exp(-2\sigma(t-s))ds \qquad (4.64)$$

Combining (4.61), (4.62) and (4.64) we obtain for all $\sigma \in [0, \lambda_1)$, $t > 0$:

$$\|u_4[t]\|_r^2 \leq \frac{1}{\lambda_1 - \sigma} \int_0^t \exp(-(\lambda_1 - \sigma)(t-s))\|f[s]\|_r^2 \exp(-2\sigma(t-s))ds$$

which directly implies

$$\|u_4[t]\|_r^2 \leq \frac{1}{(\lambda_1 - \sigma)^2} \sup_{0 < s < t} \left( \|f[s]\|_r^2 \exp(-2\sigma(t-s)) \right) \qquad (4.65)$$

Therefore, we conclude that the following inequality holds for all $\sigma \in [0, \lambda_1)$, $t > 0$:

$$\|u_4[t]\|_r \leq \frac{1}{\lambda_1 - \sigma} \sup_{0 < s < t} \left( \|f[s]\|_r \exp(-\sigma(t-s)) \right) \qquad (4.66)$$

Equation (4.44) implies that

$$\|u[t]\|_r \leq \|u_1[t]\|_r + \|u_2[t]\|_r + \|u_3[t]\|_r + \|u_4[t]\|_r, \text{ for all } t \geq 0 \qquad (4.67)$$

Using (4.67), (4.47), (4.53), (4.59) and (4.66) we get estimate (2.8). The proof is complete. ◁

**Proof of Proposition 2.6:** It suffices to show that there exists a constant $\mu > 0$

$$\int_0^1 g(x)f(x)(Af)(x)dx \geq \mu \int_0^1 f^2(x)dx, \text{ for all } f \in D \qquad (4.68)$$

Since $p(x) \equiv p$, $r(x) \equiv 1$, we obtain from (2.1) for all $f \in D$:

$$\int_0^1 g(x)f(x)(Af)(x)dx = pg(1)\left( \frac{g'(1)}{2g(1)} f(1) - f'(1) \right) f(1) + pg(0)\left( f'(0) - \frac{g'(0)}{2g(0)} f(0) \right) f(0)$$
$$+ p \int_0^1 g(x)(f'(x))^2 dx + \int_0^1 \left( g(x)q(x) - \frac{p}{2} g''(x) \right) f^2(x)dx \qquad (4.69)$$

Since $f(z) = f(0) + \int_0^z f'(s)ds$ for all $x \in [0,1]$ and since $\varepsilon_0 > 0$, we get for all $x \in [0,1]$:

$$f^2(x) \leq (1 + \varepsilon_0)f^2(0) + (1 + \varepsilon_0^{-1})\left( \int_0^x f'(s)ds \right)^2 \qquad (4.70)$$

Using the Cauchy-Schwarz inequality, which gives

$$\left(\int_0^x f'(s)ds\right)^2 \le \left(\int_0^x \frac{ds}{g(s)}\right)\left(\int_0^x g(s)(f'(s))^2 ds\right) \le \left(\int_0^x \frac{ds}{g(s)}\right)\left(\int_0^1 g(s)(f'(s))^2 ds\right), \text{ for all } x \in [0,1]$$

we get from (4.70) for all $x \in [0,1]$:

$$f^2(x) \le (1+\varepsilon_0)f^2(0) + (1+\varepsilon_0^{-1})\left(\int_0^x \frac{ds}{g(s)}\right)\left(\int_0^1 g(s)(f'(s))^2 ds\right) \tag{4.71}$$

Inequality (4.71) implies the following inequality:

$$\int_0^1 f^2(x)dx \le (1+\varepsilon_0)f^2(0) + \left((1+\varepsilon_0^{-1})\int_0^1\int_0^x \frac{ds}{g(s)}dx\right)\left(\int_0^1 g(s)(f'(s))^2 ds\right) \tag{4.72}$$

Moreover, since $f(x) = f(1) - \int_x^1 f'(s)ds$ for all $x \in [0,1]$ and since $\varepsilon_1 > 0$, we get for all $x \in [0,1]$:

$$f^2(x) \le (1+\varepsilon_1)f^2(1) + (1+\varepsilon_1^{-1})\left(\int_x^1 f'(s)ds\right)^2 \tag{4.73}$$

Using again the Cauchy-Schwarz inequality, which gives

$$\left(\int_x^1 f'(s)ds\right)^2 \le \left(\int_x^1 \frac{ds}{g(s)}\right)\left(\int_x^1 g(s)(f'(s))^2 ds\right) \le \left(\int_x^1 \frac{ds}{g(s)}\right)\left(\int_0^1 g(s)(f'(s))^2 ds\right), \text{ for all } x \in [0,1]$$

we get from (4.73) for all $x \in [0,1]$:

$$f^2(x) \le (1+\varepsilon_1)f^2(1) + (1+\varepsilon_1^{-1})\left(\int_x^1 \frac{ds}{g(s)}\right)\left(\int_0^1 g(s)(f'(s))^2 ds\right) \tag{4.74}$$

Inequality (4.74) implies the following inequality:

$$\int_0^1 f^2(x)dx \le (1+\varepsilon_1)f^2(1) + \left((1+\varepsilon_1^{-1})\int_0^1\int_x^1 \frac{ds}{g(s)}dx\right)\left(\int_0^1 g(s)(f'(s))^2 ds\right) \tag{4.75}$$

Multiplying (4.72) by $\lambda$ and (4.75) by $1-\lambda$ and adding, we get from (2.13):

$$\int_0^1 f^2(x)dx \le \lambda(1+\varepsilon_0)f^2(0) + (1-\lambda)(1+\varepsilon_1)f^2(1) + R\left(\int_0^1 g(s)(f'(s))^2 ds\right) \tag{4.76}$$

Combining (4.69) and (4.76), we get for all $f \in D$:

$$\int_0^1 g(x)f(x)(Af)(x)dx \ge pg(1)\left(\frac{g'(1)}{2g(1)}f(1) - f'(1) - \frac{(1-\lambda)(1+\varepsilon_1)}{Rg(1)}f(1)\right)f(1)$$
$$+ pg(0)\left(f'(0) - \frac{g'(0)}{2g(0)}f(0) - \frac{\lambda(1+\varepsilon_0)}{Rg(0)}f(0)\right)f(0) \tag{4.77}$$
$$+ \int_0^1\left(g(x)q(x) - \frac{p}{2}g''(x) + pR^{-1}\right)f^2(x)dx$$

Using (2.2), (2.11), (2.14), (2.15), it follows that for all $f \in D$ the following inequalities hold:

$$\left(\frac{g'(1)}{2g(1)}f(1) - f'(1) - \frac{(1-\lambda)(1+\varepsilon_1)}{Rg(1)}f(1)\right)f(1) \ge 0$$

$$\left(f'(0) - \frac{g'(0)}{2g(0)}f(0) - \frac{\lambda(1+\varepsilon_0)}{Rg(0)}f(0)\right)f(0) \ge 0$$

The above inequalities in conjunction with (4.77) imply that (4.68) holds for all $f \in D$ with $\mu := \min_{0 \leq x \leq 1}\left(g(x)q(x) - \frac{p}{2}g''(x) + pR^{-1}\right)$. Inequality (2.12) implies that $\mu > 0$. Inequality (4.68) with $f = \phi_1$ gives $\lambda_1 \geq \frac{\min\{2g(x)q(x) + 2pR^{-1} - pg''(x) : x \in [0,1]\}}{2\max\{g(x) : x \in [0,1]\}}$. The proof is complete. ◁

Since the proof of Theorem 2.7 requires Proposition 2.11, we first prove Proposition 2.11.

**Proof of Proposition 2.11:** Since the SL operator $A': D' \to L^2(0,1)$ defined by (21), (22) with $q(x) \equiv q$, satisfies assumption (H) and since $f' \in GD$ and $f(t,1)$, $f(t,0)$ are of class $GB$, it follows that that there exists an increasing sequence of times $\{\tau_i \geq 0, i = 0,1,2,...\}$ with $\tau_0 = 0$, $\lim_{i \to +\infty}(\tau_i) = +\infty$ and a unique function $v: \Re_+ \times [0,1] \to \Re$ for which the mapping $\Re_+ \ni t \to v[t] \in L^2(0,1)$ is continuous, with $v \in C^1(I \times [0,1])$ satisfying $v[t] \in C^2([0,1])$ for all $t > 0$, $\lim_{t \to \tau_i^-}(v(t,x)) = v(\tau_i, x)$, $\lim_{t \to \tau_i^-}\left(\frac{\partial v}{\partial t}(t,x)\right) = -(A'v[\tau_i])(x) + \lim_{t \to \tau_i^-}(f'(t,x))$, $\lim_{t \to \tau_i^-}\left(\frac{\partial v}{\partial x}(t,x)\right) = \frac{\partial v}{\partial x}(\tau_i, x)$, $v(0,x) = v_0(x) = u_0'(x)$ for all $x \in [0,1]$, and

$$\frac{\partial v}{\partial t}(t,x) = p\frac{\partial^2 v}{\partial x^2}(t,x) - qv(t,x) + f'(t,x), \text{ for all } (t,x) \in I \times (0,1) \quad (4.78)$$

$$\frac{\partial v}{\partial x}(t,0) + p^{-1}f(t,0) = \frac{\partial v}{\partial x}(t,1) + p^{-1}f(t,1) = 0, \text{ for all } t \in I \quad (4.79)$$

where $I = \Re_+ \setminus \{\tau_i \geq 0, i = 0,1,2,...\}$. It follows from (4.78) and (4.79) that

$$\frac{d}{dt}\int_0^1 y(t,x)dx = -q\int_0^1 y(t,x)dx, \text{ for all } t \in I \quad (4.80)$$

Since $v_0 = u_0'$ with $u_0(0) = u_0(1) = 0$ and since the mapping $t \to \int_0^1 v(t,z)dz$ is continuous, it follows that (2.26) holds. Define:

$$\tilde{u}(t,x) = \int_0^x v(t,s)ds, \text{ for } (t,x) \in \Re_+ \times [0,1] \quad (4.81)$$

It is straightforward (using (2.26), (4.78), (4.79) and definition (4.81)) to verify that

$$\frac{\partial \tilde{u}}{\partial t}(t,x) = p\frac{\partial^2 \tilde{u}}{\partial x^2}(t,x) - q\tilde{u}(t,x) + f(t,x), \text{ for all } (t,x) \in I \times (0,1) \quad (4.82)$$

$$\tilde{u}(t,0) = \tilde{u}(t,1) = 0, \text{ for all } t \geq 0 \quad (4.83)$$

Moreover, since $v_0 = u_0'$, we obtain from definition (4.81) that $\tilde{u}(0,x) = u_0(x)$ for all $x \in [0,1]$. Uniqueness of the evolution problem (4.78), (4.79) implies that $\tilde{u} \equiv u$. Equation (2.25) is a direct consequence of definition (4.81) and the fact that $\tilde{u} \equiv u$. The proof is complete. ◁

**Proof of Theorem 2.7:** The SL operator $A: D \to L^2(0,1)$ defined by (2.1), (2.2) with $p(x) \equiv p$, $r(x) \equiv 1$, $a_2 = b_2 = 0$, $a_1 = b_1 = 1$, $q(x) \equiv q$, has eigenvalues $\lambda_n = pn^2\pi^2 + q$ for $n \geq 1$ and eigenfunctions $\phi_n(x) = \sqrt{2}\sin(n\pi x)$ for $n \geq 1$. Therefore, $A$ satisfies assumption (H) and is ES since $\lambda_1 = p\pi^2 + q > 0$. Theorem 2.4 implies that for every $u_0 \in L^2(0,1)$, $f \in GD$, the unique solution $u: \Re_+ \times [0,1] \to \Re$ of the evolution equation (2.4) with (2.5), $d_0(t) = d_1(t) \equiv 0$ and initial condition $u_0 \in L^2(0,1)$ satisfies estimate (2.16) for all $\sigma \in [0, p\pi^2 + q)$ and $t > 0$. Moreover, Theorem 2.4, definition (2.18) and direct computation of the functions $\tilde{u}, \bar{u}$ involved in (2.9), (2.10) gives

$$\sqrt{\sum_{n=1}^{\infty} \frac{n^2 \pi^2}{\left(pn^2\pi^2 + q\right)^2}} = p^{-1} h(q/p). \tag{4.84}$$

The SL operator $A': D' \to L^2(0,1)$ defined by $(A'f)(x) = -p f''(x) + q f(x)$, for all $f \in D'$ and $x \in (0,1)$ with $D' \subseteq H^2(0,1)$ being the set of all functions $f:[0,1] \to \Re$ for which $f'(0) = f'(1) = 0$, has eigenvalues $\lambda_0 = q$, $\lambda_n = pn^2\pi^2 + q$ for $n \geq 1$ and eigenfunctions $\psi_0(x) \equiv 1$, $\psi_n(x) = \sqrt{2} \cos(n\pi x)$ for $n \geq 1$. Therefore, $A'$ satisfies assumption (H) but it is not necessarily an ES operator. Consequently, we cannot use Theorem 2.4 to system (2.23), (2.24). However, we are in a position to prove (exactly as in the proof of Theorem 2.4) the following claim. Its proof is omitted, since it is identical to that in the proof of Theorem 2.4 and the only additional thing is the use of (2.26).

**Claim:** For every solution $v: \Re_+ \times [0,1] \to \Re$ of the evolution equation (2.23) with (2.24), initial condition $v_0 = u'_0$, where $u_0 \in H^1(0,1)$ with $u_0(0) = u_0(1) = 0$, the following equations hold for all $t \geq 0$:

$$\begin{aligned} c_n(t) &= \exp\left(-\left(pn^2\pi^2 + q\right)t\right) c_n(0) - (-1)^n \sqrt{2} \int_0^t \exp\left(-\left(pn^2\pi^2 + q\right)(t-s)\right) f(s,1) ds \\ &+ \sqrt{2} \int_0^t \exp\left(-\left(pn^2\pi^2 + q\right)(t-s)\right) f(s,0) ds + \sqrt{2} \int_0^t \exp\left(-\left(pn^2\pi^2 + q\right)(t-s)\right) \left(\int_0^1 f'(s,z) \cos(n\pi z) dz\right) ds \end{aligned}$$
, for $n = 1, 2, \ldots$ (4.85)

$$c_0(t) \equiv 0 \tag{4.86}$$

where

$$c_0(t) := \int_0^1 y(t,x) dx, \quad c_n(t) := \sqrt{2} \int_0^1 y(t,x) \cos(n\pi x) dx, \text{ for } n = 1, 2, \ldots \tag{4.87}$$

Next, we recognize that the solution $v: \Re_+ \times [0,1] \to \Re$ of the evolution equation (2.23) with (2.24), initial condition $v_0 = u'_0$, where $u_0 \in H^1(0,1)$ with $u_0(0) = u_0(1) = 0$ satisfies the following equation:

$$v[t] = v_1[t] + v_2[t], \text{ for all } t \geq 0 \tag{4.88}$$

where
- $v_1: \Re_+ \times [0,1] \to \Re$ is the solution of the evolution equation (2.23) with (2.24), initial condition $v_0 = u'_0$, where $u_0 \in H^1(0,1)$ with $u_0(0) = u_0(1) = 0$ and corresponding to input $f \equiv 0$,
- $v_2: \Re_+ \times [0,1] \to \Re$ is the solution of the evolution equation (2.23) with (2.24), initial condition $v_0 \equiv 0$ and corresponding to input $f \in GD$.

Next, we estimate each component of $v: \Re_+ \times [0,1] \to \Re$ separately.

1) Estimate for $v_1: \Re_+ \times [0,1] \to \Re$.
We obtain from (4.85) for all $t \geq 0$ and $n = 1, 2, \ldots$:

$$c_n(t) = \exp\left(-\left(pn^2\pi^2 + q\right)t\right) c_n(0) \tag{4.89}$$

where $c_n(t) := \sqrt{2} \int_0^1 y(t,x) \cos(n\pi) dx$. Since the eigenfunctions $\{\psi_n\}_{n=0}^{\infty}$ of the SL operator $A': D' \to L^2(0,1)$ defined by $(A'f)(x) = -p f''(x) + q f(x)$, for all $f \in D'$ and $x \in (0,1)$ with $D' \subseteq H^2(0,1)$ being the set of all functions $f:[0,1] \to \Re$ for which $f'(0) = f'(1) = 0$, form an orthonormal basis of $L^2(0,1)$, it follows that Parseval's identity holds, i.e.,

$$\|v_1[t]\|^2 = c_0^2(t) + \sum_{n=1}^{\infty} c_n^2(t), \text{ for all } t \geq 0 \tag{4.90}$$

Using (4.86), (4.89) and the fact that $v_0 = u'_0$, we get for all $t \geq 0$:

$$\|v_1[t]\| \leq \exp\left(-\left(pn^2\pi^2 + q\right)t\right) \|u'_0\| \tag{4.91}$$

2) Estimate for $v_2 : \mathfrak{R}_+ \times [0,1] \to \mathfrak{R}$.

Using integration by parts, we obtain from (4.85) for all $t \geq 0$ and $n = 1, 2, \ldots$:

$$c_n(t) = n\pi \int_0^t \exp\left(-(pn^2\pi^2 + q)(t-s)\right) f_n(s) ds \tag{4.92}$$

where $c_n(t) := \sqrt{2} \int_0^1 y(t,x) \cos(n\pi x) dx$, $f_n(t) := \sqrt{2} \int_0^1 f(t,x) \sin(n\pi x) dx$. Since $\{\psi_n\}_{n=0}^\infty$, $\{\phi_n\}_{n=1}^\infty$ are orthonormal bases of $L^2(0,1)$, it follows that Parseval's identity holds, i.e.,

$$\|v_2[t]\|^2 = c_0^2(t) + \sum_{n=1}^\infty c_n^2(t), \text{ for all } t \geq 0 \tag{4.93}$$

$$\|f[t]\|^2 = \sum_{n=1}^\infty f_n^2(t), \text{ for all } t \geq 0 \tag{4.94}$$

It follows from (4.92) that the following inequality holds for all $\sigma \in [0, p\pi^2 + q)$, $t \geq 0$ and $n = 1, 2, \ldots$:

$$|c_n(t)| \leq n\pi \int_0^t \exp\left(-(pn^2\pi^2 + q)(t-s)\right) |f_n(s)| ds \leq n\pi \int_0^t \exp\left(-(pn^2\pi^2 + q - \sigma)(t-s)\right) |f_n(s)| \exp(-\sigma(t-s)) ds \tag{4.95}$$

Using the Cauchy-Schwarz inequality and (4.95) we obtain for all $\sigma \in [0, p\pi^2 + q)$, $t \geq 0$ and $n = 1, 2, \ldots$:

$$|c_n(t)| \leq \frac{n\pi}{\sqrt{pn^2\pi^2 + q - \sigma}} \left( \int_0^t \exp\left(-(pn^2\pi^2 + q - \sigma)(t-s)\right) |f_n(s)|^2 \exp(-2\sigma(t-s)) ds \right)^{1/2}$$

which directly implies

$$|c_n(t)|^2 \leq \frac{n^2\pi^2}{pn^2\pi^2 + q - \sigma} \int_0^t \exp\left(-(pn^2\pi^2 + q - \sigma)(t-s)\right) |f_n(s)|^2 \exp(-2\sigma(t-s)) ds \tag{4.96}$$

Equation (4.94) implies that $\|f[t]\|^2 \geq f_n^2(t)$ for all $t \geq 0$ and $n = 1, 2, \ldots$. Combining the previous inequality with (4.96) we obtain for all $\sigma \in [0, p\pi^2 + q)$, $t \geq 0$ and $n = 1, 2, \ldots$:

$$|c_n(t)|^2 \leq \frac{n^2\pi^2}{pn^2\pi^2 + q - \sigma} \int_0^t \exp\left(-(pn^2\pi^2 + q - \sigma)(t-s)\right) \|f[s]\|^2 \exp(-2\sigma(t-s)) ds \tag{4.97}$$

from which we get for all $\sigma \in [0, p\pi^2 + q)$, $t > 0$ and $n = 1, 2, \ldots$:

$$|c_n(t)|^2 \leq \frac{n^2\pi^2}{(pn^2\pi^2 + q - \sigma)^2} \sup_{0 < s < t} \left( \|f[s]\|^2 \exp(-2\sigma(t-s)) \right)$$

The above inequality in conjunction with (4.86) and (4.93) gives for all $\sigma \in [0, p\pi^2 + q)$, $t > 0$:

$$\|v_2[t]\|^2 \leq \left( \sum_{n=1}^\infty \frac{n^2\pi^2}{(pn^2\pi^2 + q - \sigma)^2} \right) \sup_{0 < s < t} \left( \|f[s]\|^2 \exp(-2\sigma(t-s)) \right)$$

$$\leq \left( \sum_{n=1}^\infty \left( \frac{pn^2\pi^2 + q}{pn^2\pi^2 + q - \sigma} \right)^2 \frac{n^2\pi^2}{(pn^2\pi^2 + q)^2} \right) \sup_{0 < s < t} \left( \|f[s]\|^2 \exp(-2\sigma(t-s)) \right)$$

$$\leq \left( \frac{p\pi^2 + q}{p\pi^2 + q - \sigma} \right)^2 \left( \sum_{n=1}^\infty \frac{n^2\pi^2}{(pn^2\pi^2 + q)^2} \right) \sup_{0 < s < t} \left( \|f[s]\|^2 \exp(-2\sigma(t-s)) \right)$$

which combined with (4.84) directly implies that

$$\|v_2[t]\| \leq \left( \frac{p\pi^2 + q}{p\pi^2 + q - \sigma} \right) p^{-1} h(p/q) \sup_{0 < s < t} \left( \|f[s]\| \exp(-\sigma(t-s)) \right) \tag{4.98}$$

Equation (4.88) implies that

$$\|v[t]\| \leq \|v_1[t]\| + \|v_2[t]\|, \text{ for all } t \geq 0 \tag{4.99}$$

Using (2.25), (4.99), (4.91) and (4.98) we get estimate (2.17). The proof is complete. ◁

The proof of Theorem 2.9 requires Proposition 2.12 and the following technical proposition, which we prove next.

**Proposition 4.4:** *Let* $T > 0$, $p > 0$, $a_1, a_2, b_1, b_2$ *be real constants with* $|a_1| + |a_2| > 0$, $|b_1| + |b_2| > 0$ *and let* $q \in C^0([0,1])$, $\beta \in C^0([0,1]^2)$ *be given functions. Let* $u : [0,T] \times [0,1] \to \Re$ *be a function for which the mapping* $[0,T] \ni t \to u[t] \in L^2(0,1)$ *is continuous and for which there exists a finite set* $\{\tau_i \in (0,T), i = 0,1,2,...,N\}$ *such that* $u \in C^1(I \times [0,1])$, *where* $I = (0,T) \setminus \{\tau_i \in (0,T), i = 0,1,2,...,N\}$. *Finally, suppose that* $u[0] \equiv 0$, $u[t] \in C^2([0,1])$ *for all* $t \in I$ *and that*

$$\frac{\partial u}{\partial t}(t,x) = p \frac{\partial^2 u}{\partial x^2}(t,x) - q(x) u(t,x) + \int_0^x \beta(x,s) u(t,s) ds, \text{ for } x \in (0,1) \text{ and } t \in I \tag{4.100}$$

$$b_2 \frac{\partial u}{\partial x}(t,0) + b_1 u(t,0) = a_2 \frac{\partial u}{\partial x}(t,1) + a_1 u(t,1) = 0, \text{ for } t \in I \tag{4.101}$$

*Then* $u[t] \equiv 0$ *for all* $t \in [0,T]$.

**Proof:** Let $k \in C^2([0,1]; (0,+\infty))$ be a positive function that satisfies

$$-\frac{k'(1)}{2k(1)} \leq \begin{cases} a_1/a_2 & \text{if } a_2 \neq 0 \\ +\infty & \text{if } a_2 = 0 \end{cases}, \quad -\frac{k'(0)}{2k(0)} \geq \begin{cases} b_1/b_2 & \text{if } b_2 \neq 0 \\ -\infty & \text{if } b_2 = 0 \end{cases} \tag{4.102}$$

Consider the mapping

$$V(t) = \frac{1}{2} \int_0^1 k(x) u^2(t,x) dx, \text{ for } t \in [0,T] \tag{4.103}$$

Notice that continuity of the mapping $[0,T] \ni t \to u[t] \in L^2(0,1)$ implies continuity of the mapping $[0,T] \ni t \to V(t) \in \Re$. The fact that $u \in C^1(I \times [0,1])$ implies that $V(t)$ is $C^1$ on $I$ and satisfies

$$\dot{V}(t) = \int_0^1 k(x) u(t,x) \frac{\partial u}{\partial t}(t,x) dx, \text{ for } t \in I \tag{4.104}$$

Using (4.104), (4.100), (4.101), (4.103) the Cauchy-Schwarz inequality, repeated integration by parts and (4.102), we obtain for all $t \in I$:

$$\dot{V}(t) = p \int_0^1 k(x) u(t,x) \frac{\partial^2 u}{\partial x^2}(t,x) dx - \int_0^1 k(x) q(x) u^2(t,x) dx + \int_0^1 k(x) u(t,x) \left( \int_0^x \beta(x,s) u(t,s) ds \right) dz$$

$$\leq \max_{0 \leq x \leq 1}(-q(x)) \int_0^1 k(x) u^2(t,x) dx + pk(1) \left( \frac{\partial u}{\partial x}(t,1) - \frac{k'(1)}{2k(1)} u(t,1) \right) u(t,1)$$

$$- pk(0) \left( \frac{\partial u}{\partial x}(t,0) - \frac{k'(0)}{2k(0)} u(t,0) \right) u(t,0) + \frac{p}{2} \int_0^1 k''(x) u^2(t,x) dx$$

$$- p \int_0^1 k(x) \left( \frac{\partial u}{\partial x}(t,x) \right)^2 dx + \left( \int_0^1 k(x) u^2(t,x) dx \right)^{1/2} \left( \int_0^1 k(x) \left( \int_0^x k(s) u^2(t,s) ds \right) \left( \int_0^x \frac{\beta^2(x,s)}{k(s)} ds \right) dx \right)^{1/2}$$

$$\leq \left( \max_{0 \leq x \leq 1}(-q(x)) + \max_{0 \leq x \leq 1}\left( \frac{pk''(x)}{2k(x)} \right) + \left( \int_0^1 k(x) \left( \int_0^x \frac{\beta^2(x,s)}{k(s)} ds \right) dx \right)^{1/2} \right) \int_0^1 k(x) u^2(t,x) dx$$

$$+ pk(1) \left( \frac{\partial u}{\partial x}(t,1) - \frac{k'(1)}{2k(1)} u(t,1) \right) u(t,1) - pk(0) \left( \frac{\partial u}{\partial x}(t,0) - \frac{k'(0)}{2k(0)} u(t,0) \right) u(t,0)$$

$$\leq 2 \left( \max_{0 \leq x \leq 1}(-q(x)) + \max_{0 \leq x \leq 1}\left( \frac{pk''(x)}{2k(x)} \right) + \left( \int_0^1 k(x) \left( \int_0^x \frac{\beta^2(x,s)}{k(s)} ds \right) dx \right)^{1/2} \right) V(t)$$

It follows that there exists a constant $L > 0$ such that $\dot{V}(t) \leq LV(t)$ for all $t \in I$. Continuity of the mapping $[0,T] \ni t \to V(t) \in \Re$ in conjunction with the differential inequality $\dot{V}(t) \leq LV(t)$ for $t \in I$ and the fact that $u[0] \equiv 0$ implies that

$$V(t) \leq L\int_0^t V(s)ds \text{ for all } t \in [0,T] \tag{4.105}$$

Using Gronwall's inequality Lemma and (4.105), we get $V(t) \leq 0$, for all $t \in [0,T]$, which combined with definition (4.103) gives $u[t] \equiv 0$ for all $t \in [0,T]$. The proof is complete. ◁

**Proof of Proposition 2.12:** Since $f \in GD$, $d_1 \in GB$ it follows from Theorem 2.2 that for every $x_0 \in H^1(0,1)$ there exist an increasing sequence of times $\{\tau_i \geq 0, i = 0,1,2,...\}$ with $\tau_0 = 0$, $\lim_{i \to +\infty}(\tau_i) = +\infty$ and a unique function $u : \Re_+ \times [0,1] \to \Re$ for which the mapping $\Re_+ \ni t \to u[t] \in L^2(0,1)$ is continuous, with $u \in C^1(I \times [0,1])$ satisfying $u[t] \in C^2([0,1])$ for all $t > 0$, $\lim_{t \to \tau_i^-}(u(t,x)) = u(\tau_i, x)$, $\lim_{t \to \tau_i^-}\left(\frac{\partial u}{\partial t}(t,x)\right) = -(Au[\tau_i])(x) + \lim_{t \to \tau_i^-}(f(t,x))$, $\lim_{t \to \tau_i^-}\left(\frac{\partial u}{\partial x}(t,x)\right) = \frac{\partial u}{\partial x}(\tau_i, x)$, $u(0,x) = u_0(x)$ for all $x \in [0,1]$, and equations (2.6), (2.7) hold, where $I = \Re_+ \setminus \{\tau_i \geq 0, i = 0,1,2,...\}$.

Since $f \in GD$, $f' \in GD$ and $q \in C^2([0,1])$, it follows from Definition 2.1 that the mapping $\tilde{f} : (0,+\infty) \times [0,1] \to \Re$ defined by

$$\tilde{f}(t,x) := -\left(q'(x) - 2a_1^2 px\right)u(t,x) + f'(t,x) + a_1 x f(t,x), \text{ for } (t,x) \in \Re_+ \times [0,1] \tag{4.106}$$

is of class $GD$.

Since $\tilde{f} \in GD$, $d_1 \in GB$ and since the SL operator $A' : D' \to L^2(0,1)$ defined by (2.19), (2.20), satisfies assumption (H), it follows from Theorem 2.2 that there exists an increasing sequence of times $\{\tau'_i \geq 0, i = 0,1,2,...\}$ with $\tau'_0 = 0$, $\lim_{i \to +\infty}(\tau'_i) = +\infty$ (possibly different from $\{\tau_i \geq 0, i = 0,1,2,...\}$) such that we can construct a (unique) function $v : \Re_+ \times [0,1] \to \Re$ for which the mapping $\Re_+ \ni t \to v[t] \in L^2(0,1)$ is continuous, with $v \in C^1(I' \times [0,1])$ satisfying $v[t] \in C^2([0,1])$ for all $t > 0$, $v(0,x) = v_0(x) = u'_0(x) + a_1 x u_0(x)$ for all $x \in [0,1]$, and

$$\frac{\partial v}{\partial t}(t,x) = p\frac{\partial^2 v}{\partial x^2}(t,x) - (q(x) + 2a_1 p)v(t,x) + \tilde{f}(t,x), \text{ for all } (t,x) \in I' \times (0,1) \tag{4.107}$$

$$\frac{\partial v}{\partial x}(t,0) + p^{-1}f(t,0) = v(t,1) - d_1(t) = 0, \text{ for all } t \in I' \tag{4.108}$$

where $I' = \Re_+ \setminus \{\tau'_i \geq 0, i = 0,1,2,...\}$.

Define

$$\tilde{u}(t,x) := \int_0^x \exp\left(-\frac{a_1}{2}(x^2 - s^2)\right)v(t,s)ds, \text{ for } (t,x) \in \Re_+ \times [0,1] \tag{4.109}$$

Using repeated integration by parts, (4.107), (4.108), (4.109) and the facts that $v(0,x) = v_0(x) = u'_0(x) + a_1 x u_0(x)$, $u_0(0) = 0$, we show that:

$$\frac{\partial \tilde{u}}{\partial t}(t,x) = p\frac{\partial^2 \tilde{u}}{\partial x^2}(t,x) - q(x)\tilde{u}(t,x) + \int_0^x \exp\left(-\frac{a_1}{2}(x^2 - s^2)\right)\left(q'(s) - 2a_1^2 ps\right)(\tilde{u}(t,s) - u(t,s))ds + f(t,x),$$

for all $(t,x) \in I' \times (0,1)$ \hfill (4.110)

$$\tilde{u}(t,0) = \frac{\partial \tilde{u}}{\partial x}(t,1) + a_1\tilde{u}(t,1) - d_1(t) = 0, \text{ for all } t \in I' \tag{4.111}$$

$$v(t,z) = \frac{\partial \tilde{u}}{\partial x}(t,x) + a_1 x \tilde{u}(t,x), \text{ for all } (t,x) \in (0,+\infty) \times [0,1] \quad (4.112)$$

$$\tilde{u}(0,x) = u_0(x), \text{ for all } x \in [0,1]. \quad (4.113)$$

Notice that definition (4.109) and continuity of the mapping $\Re_+ \ni t \to v[t] \in L^2(0,1)$ implies that the mapping $\Re_+ \ni t \to \tilde{u}[t] \in L^2(0,1)$ is continuous. Applying Proposition 4.4 to the function $\bar{u} = \tilde{u} - u$ and using (2.6), (2.7), (4.110), (4.111) and (4.113) we get $u[t] = \tilde{u}[t]$ for $t \geq 0$. Equation (2.29) is a direct consequence of (4.112), the fact that $u[t] = \tilde{u}[t]$ for $t \geq 0$ and the fact that $v(0,x) = v_0(x) = u_0'(x) + a_1 x u_0(x)$. The proof is complete. ◁

We are now in a position to prove Theorem 2.9.

**Proof of Theorem 2.9:** Theorem 2.4 implies that for every $f \in GD$, $d_1 \in GB$, $u_0 \in L_r^2(0,1)$, the unique solution $u : \Re_+ \times [0,1] \to \Re$ of the evolution equation (2.4) with (2.5) and initial condition $u_0 \in H^1(0,1)$ satisfies the following estimate for all $\sigma \in [0, \lambda_1)$ and $t > 0$:

$$\|u[t]\| \leq \exp(-\sigma t)\|u_0\| + \frac{\lambda_1}{\lambda_1 - \sigma} C_1 \sup_{0 < s < t}(|d_1(s)|\exp(-\sigma(t-s))) + \frac{1}{\lambda_1 - \sigma} \sup_{0 < s < t}(\|f[s]\|\exp(-\sigma(t-s))) \quad (4.114)$$

where

$$C_0 := \|\tilde{u}\|, \ C_1 := \frac{1}{\sqrt{a_1^2 + 1}} \|\bar{u}\| \quad (4.115)$$

$\tilde{u} \in C^2([0,1])$ is the unique solution of the boundary value problem $p\tilde{u}''(x) - q(x)\tilde{u}(x) = 0$ for $x \in [0,1]$ with $\tilde{u}(0) = 1$, $a_1\tilde{u}(1) + \tilde{u}'(1) = 0$ and $\bar{u} \in C^2([0,1])$ is the unique solution of the boundary value problem $p\bar{u}''(x) - q(x)\bar{u}(x) = 0$ for $x \in [0,1]$ with $\bar{u}(0) = 0$ and $a_1\bar{u}(1) + \bar{u}'(1) = \sqrt{a_1^2 + 1}$.

Let $0 < \mu_1 < \mu_2 < \ldots < \mu_n < \ldots$ with $\lim_{n \to \infty}(\mu_n) = +\infty$ be the eigenvalues of the SL operator $A'$. Then it follows from Proposition 2.12 and Theorem 2.4 that the unique solution $v : \Re_+ \times [0,1] \to \Re$ of the evolution equation (2.27) with (2.28) and initial condition $v_0(x) = u_0'(x) + a_1 x u_0(x)$ for $x \in [0,1]$ satisfies the following estimate for all $\sigma \in [0, \mu_1)$ and $t > 0$:

$$\|v[t]\| \leq \exp(-\sigma t)\|v_0\| + \frac{\mu_1}{\mu_1 - \sigma} \tilde{C}_0 p^{-1} \sup_{0 < s < t}(|f(s,0)|\exp(-\sigma(t-s)))$$
$$+ \frac{\mu_1}{\mu_1 - \sigma} \tilde{C}_1 \sup_{0 < s < t}(|d_1(s)|\exp(-\sigma(t-s))) + \frac{1}{\mu_1 - \sigma} \sup_{0 < s < t}(\|\tilde{f}[s]\|\exp(-\sigma(t-s))) \quad (4.116)$$

where the mapping $\tilde{f} : \Re_+ \times [0,1] \to \Re$ is defined by (4.106) and

$$\tilde{C}_0 := \|\tilde{v}\|, \ \tilde{C}_1 := \|\bar{v}\| \quad (4.117)$$

$\tilde{v} \in C^2([0,1])$ is the unique solution of the boundary value problem $p\tilde{v}''(x) - (q(x) + 2a_1 p)\tilde{v}(x) = 0$ for $x \in [0,1]$ with $\tilde{v}'(0) = 1$, $\tilde{v}(1) = 0$ and $\bar{v} \in C^2([0,1])$ is the unique solution of the boundary value problem $p\bar{v}''(x) - (q(x) + 2a_1 p)\bar{v}(x) = 0$ for $x \in [0,1]$ with $\bar{v}'(0) = 0$ and $\bar{v}(1) = 1$. Using (2.29), (4.106) and the fact that $v_0(x) = u_0'(x) + a_1 x u_0(x)$ for $x \in [0,1]$, we obtain the following inequalities:

$$\|\tilde{f}[t]\| \leq \max_{0 \leq x \leq 1}(|q'(x) - 2a_1^2 p x|)\|u[t]\| + \|f'[t]\| + |a_1|\|f[t]\|, \text{ for all } t > 0 \quad (4.118)$$

$$\|u'[t]\| \leq \|v[t]\| + |a_1|\|u[t]\|, \text{ for all } t > 0 \quad (4.119)$$

$$\|v_0\| \leq \|u_0'\| + |a_1|\|u_0\| \quad (4.120)$$

It follows from (4.114), (4.116), (4.118), (4.119) and (4.120) that the following inequalities hold for all $\sigma \in [0, \min(\mu_1, \lambda_1))$ and $t > 0$:

$$\|u'[t]\| \leq \exp(-\sigma t)\|u'_0\| + \exp(-\sigma t)\left(2|a_1| + \frac{\max_{0\leq x\leq 1}\left(|q'(x) - 2a_1^2 px|\right)}{\mu_1 - \sigma}\right)\|u_0\|$$

$$+ \frac{\mu_1}{\mu_1 - \sigma}\tilde{C}_0 p^{-1} \sup_{0<s<t}\left(|f(s,0)|\exp(-\sigma(t-s))\right)$$

$$+ \left(\frac{\mu_1}{\mu_1 - \sigma}\tilde{C}_1 + \frac{\lambda_1}{\lambda_1 - \sigma}C_1\left(\frac{\max_{0\leq x\leq 1}\left(|q'(x) - 2a_1^2 px|\right)}{\mu_1 - \sigma} + |a_1|\right)\right)\sup_{0<s<t}\left(|d_1(s)|\exp(-\sigma(t-s))\right) \quad (4.121)$$

$$+ \frac{1}{\mu_1 - \sigma}\sup_{0<s<t}\left(\|f'[s]\|\exp(-\sigma(t-s))\right)$$

$$+ \left(\frac{|a_1|}{\mu_1 - \sigma} + \frac{|a_1|}{\lambda_1 - \sigma} + \lambda_1 \frac{\max_{0\leq x\leq 1}\left(|q'(x) - 2a_1^2 px|\right)}{(\mu_1 - \sigma)(\lambda_1 - \sigma)}\right)\sup_{0<s<t}\left(\|f[s]\|\exp(-\sigma(t-s))\right)$$

Estimates (2.21), (2.22) are consequences of (4.114), (4.121). The proof is complete. ◁

## 5. Conclusions

In this work ISS estimates were derived for the solutions of 1-D linear parabolic PDEs with disturbances at both boundaries and distributed disturbances. The decay estimates were expressed in the $L^2$ and $H^1$ norms of the solution and discontinuous disturbances are allowed. The obtained estimates do not require knowledge of the eigenvalues and the eigenfunctions of the corresponding Sturm-Liouville operator and can be applied in a straightforward way for the stability analysis of parabolic PDEs with nonlocal terms.

Future work may include the study of necessary conditions for ISS in the $H^1$ norm of the state. While Theorem 2.7 provides necessary and sufficient conditions for ISS, this is not the case with Theorem 2.9, where only sufficient conditions for ISS are provided. Finally, novel results will be needed for the derivation of ISS estimates in the $H^2$ norm of the state.